\documentclass{article}
\usepackage{enumerate,amsmath,amssymb,amsthm}
\usepackage{authblk}
\usepackage{enumerate,amsmath,amssymb,amsthm}
\newtheorem{theorem}{Theorem}[section]
\newtheorem{lemma}[theorem]{Lemma}
\newtheorem{remark}[theorem]{Remark}
\newtheorem{definition}[theorem]{Definition}
\newtheorem{proposition}[theorem]{Proposition}
\newtheorem{example}[theorem]{Example}

\newtheorem{assumption}[theorem]{Assumption}
\def\bt{\begin{theorem}}
\def\et{\end{theorem}}
\def\bl{\begin{lemma}}
\def\el{\end{lemma}}
\def\br{\begin{remark}}
\def\er{\end{remark}}
\def\bd{\begin{definition}}
\def\ed{\end{definition}}
\def\bp{\begin{proposition}}
\def\ep{\end{proposition}}

\author{Sheng Wang}
\title{Approximation Theorems For Reflected Stochastic\\
Differential Equations}
\affil{School of Mathematics, Sun Yat-sen University,\\
Guangzhou, Guangdong 510275, P.R.China\\
Email: wangsh325@mail2.sysu.edu.cn}
\date{}

\begin{document}
\maketitle

\begin{abstract}
In this paper we prove a general approximation result for reflected stochastic differential equations in bounded domains satisfying conditions reorganized by Ren and Wu $\cite{rw}$. Then we show that it includes Wong-Zakai approximation, mollifier approximation, etc. 
\end{abstract}

$\textbf{Keywords}$: Mollifier Approximation; Wong-Zakai Approximation; Reflected SDEs.

\section{Introduction}
~~~~~The approximation of stochastic differential equations on a domain D has been a much studied probelm after the work of Wong and Zakai $\cite{wz}$. Wong and Zakai found the relationship between stochastic differential equations and ordinary differential equations. Their result shows that if the Brownian motion in stochastic differential equation is replaced by its linear interpolation, we get an ordinary differential equation which can approximate the corresponding Stratonovish stochastic differential equation. The approximation is the so-called Wong-Zakai approximation. This result was improved by Ikeda, Nakao and Yamato $\cite{Nakao}$ and by Ikeda and Watanabe $\cite{iw}$. They extended linear interpolation to more general Wiener transform. They can do this, since $\sigma(X^{\delta}_{s})dB^{\delta}_{s}$ converges in some way under their assumption (Assumption $\ref{a1}$ below). When it comes to the approximation of the reflected stochastic differential equation, it's necessary to mention the work of Zhang $\cite{zhangtusheng}$ and the work of Aida and Sasaki $\cite{aida}$. They established the Wong-Zakai approximation of reflected stochastic differential equations respectively, but Zhang employed an adapted version of the Wong-Zakai approximation which is helpful when people use BDG's inequality and Ito's formula.

Existence and uniqueness of reflected stochastic differential equation were proved by Lions and Sznitman $\cite{ls}$ when D is a bounded open set satisfying admissible condition. After that, this result were extended to the unbounded domain by Saisho $\cite{Saisho}$, and he removed the admissible condition. Here admissible means that D can be approximated in some sense by smooth domains.

However, the general approximation of reflected differential equation based on Assumption $\ref{a1}$ has not yet been touched. In this paper we get such an approximation
and show that it includes Wong-Zakai approximation, mollifier approximation, etc. The way we employ is similar to Zhang's. Actually, to bound $dB^{\delta}_{u}dB^{\delta}_{s}$, we must use It\^o's formula at a length of $\widetilde{\delta}$ ($=n\delta$, see definition before Lemma $\ref{l8}$), and it makes another method no longer doable. We refer to the paper of Ren and Wu $\cite{rw}$ when dealing with  $dK^{\delta}_{u}dB^{\delta}_{s}$. They have cited Evans and Stroock's work $\cite{stroock}$ in their proof, and it also brings us inspiration.

In Section 1 we give some basic definitions and state our main result under an important assumption. There are several remarks attached below that will be useful in Section 3. In Section 2, we restate the framework reorganized by Ren and Wu $\cite{rw}$, and introduce some results from others' work. Next, we provide in Section 3 some moment estimates and the proof of main result. Finally in Section 4, our examples show that the shifted approximation of Brownian motion includes Wong-Zakai approximation, mollifier approximation, etc.

Henceafter, we use $C[a,b]$ to represent the space of continuous functions in $[a,b]$, and $\mathbb{C}[a,b]$ is the $\sigma-$ algebra generated by supremum norm, and $\mathbb{C}^{n}[a,b]$ means $\mathbb{C}[a,b]\times\mathbb{C}[a,b]\cdots \mathbb{C}[a,b]$. Also $f(\ast\wedge T)$ is represented by $f^{T}(\ast)$ and $A\lesssim B$ means that there exists a $C\geq 0$ such that  $|A|\leq C|B|$. For any $f\in C[0,\infty)$, we employ notations:
\begin{equation}
\begin{aligned}
&\|f\|:=\sup _{u\in[0, \infty)}\left|f_{u}\right|,\\
&\|f\|_{C[a,b]}:=\sup _{u\in[a,b]}\left|f_{u}\right|,\\
&\|f\|_{[s, t]}:=\sup _{u, v \in[s, t]}\left|f_{u}-f_{v}\right|,\\
&\|f\|_{[s, t],\theta}:=\sup _{u, v \in[s, t],u\neq v}\frac{\left|f_{u}-f_{v}\right|}{|u-v|^{\theta}},\\
&|f|_{s}^{t}:=\sup _{\Delta} \sum_{k=1}^{N}\left|f\left(t_{k}\right)-f\left(t_{k-1}\right)\right|,
\end{aligned}
\end{equation}
where $\Delta=\left\{s=t_{0}<\cdots<t_{N}=t\right\}$ is a partition of the interval $[s, t]$. If $s\in[k\delta,~(k+1)\delta)$, we define $\underline{\underline{s}}_{\delta}:=(k-2)\delta\vee 0$, $\underline{s}_{\delta}:=(k-1)\delta\vee 0$, $s_{\delta}:=k\delta$, $\bar{s}_{\delta}:=(k+1)\delta$.

$\boldsymbol{Skorohod~problem}$. Let $D$ be a domain in $R^{d}$. There is a continuous function $h(t):[0,T]\mapsto R^{d}$ with $h(0)\in\overline{D}$. We say a pair of functions $(x(t),k(t))$ is the solution of Skorohod problem $(h,D,\mathcal{N})$, if the followng propertys hold.

(i) $x(t)=h(t)+k(t)$, for $t\in[0,T]$.

(ii) $x:[0,T]\mapsto R^{d}$ is a continuous path in $\overline{D}$ and $x(0)=h(0)$.

(iii) $k:[0,T]\mapsto R^{d}$ is a continuous function with bounded variation such that $k(0)=0$ and
$$
k(t)=\int_{0}^{t}n(s) d|k|_{0}^{s},
$$
$$
|k|_{0}^{t}=\int_{0}^{t} 1_{\partial D}(x(s)) d|k|_{0}^{s},
$$
where $n(s)\in\mathcal{N}_{x(s)}$ and $\mathcal{N}_{x}$ is the set of inward normal unit vectors at $x\in\partial D$ defined by
$$
\mathcal{N}_x=\bigcup_{r>0}\mathcal{N}_{x,r},
$$
$$
\mathcal{N}_{x,r}=\{n\in R^{d};~|n|=1,~B(x-rn,r)\cap D=\emptyset\}.
$$

$\boldsymbol{Wiener~space}$. We have a probability measure $P$ on $\left(C[0, \infty), \mathcal{C}[0, \infty)\right)$, under which the coordinate mapping process $W_{t}(w):=w(t), 0 \leq t<\infty,$ is a standard,
$r-$dimensional Brownian motion. We call $\left(C[0, \infty), \mathcal{C}[0, \infty), P\right)$ Wiener space.

$\boldsymbol{Approximation~of~Brownian~motion}$. A family $\{G^{\delta}(t, f)=(G^{\delta,1}(t, f), \ldots,~G^{\delta,r}(t, f))\}_{\delta>0}$, where $G^{\delta}~:[0,\infty)\times C[0,\infty)\mapsto~R^{r}$ is a measurable map, is called an approximation of Brownian motion, if it is a class of $r-$dimensional continuous processes defined over the Wiener space such that

(i) for every $w, t \longmapsto G^{\delta}(t, w)$ is piecewise continuously differentiable.

(ii) $G^{\delta}(t, \ast)$ is $\mathbb{C}[(t-\delta)\vee 0, t+\delta]-$measurable and is linear at $t=0$.

(iii) $G^{\delta}(t+k \delta, w)=G^{\delta}\left(t, \theta_{k\delta} w\right)+w(k\delta)$ for every $k=1,2, \dots, t \geq 0$
and $w$, where $\left(\theta_{t} w\right)(s)=w(t+s)-w(t)$.

(iv) $E[G^{\delta,i}(0, w)]=0$ for $i=1,2,\ldots,r$.

(v)
$$
\delta^{m-1}E\left[\int_{0}^{\delta}|\dot{G}^{\delta,i}(s, w)|^{2m}ds \right] \leq C \quad for \quad i=1,2, \dots, r,$$
$\forall m\geq 1$, $\delta<1$, where $\dot{G}^{\delta,i}(s, w)=\frac{d}{d s}G^{\delta,i}(s, w)$ and $C$ is a constant dependent on $m$.

(vi)
$$
E\left[\left|G^{\delta,i}(0, f(\ast))\right|^{2p}\right] \leq C_{p} \left\|E\left[\left|f(\ast)\right|^{2p}\right]\right\|_{C[0,\delta]}
$$
for  $i=1,2, \dots, r$ and $p\geq 1$, for any
$$
f(s)=g(w(t_{1}+s), \cdots,w(t_{n}+s),w),~w\in C[0,\infty),~t_{j}\in R,~n\geq0,
$$
is continuous in $[0,\delta]$, $g$ is $\mathcal{B}(R^{n})\times \mathbb{C}[0,\infty)-$measurable and $E\left[\left|f(s)\right|^{2p}\right]$ is continuous in $[0,\delta]$. At the same time, we call $B^{\delta}(t,w)$ the shifted approximation of Brownian motion, if
\begin{equation}
B^{\delta}(t,w)=
\left\{
\begin{array}{cc}
G^{\delta}(t-\delta, w), & \delta\leq t,\\
0, & 0\leq t<\delta.
\end{array}
\right.
\end{equation}

\br\label{r1}
According to (vi), let $g(x)=x$, $n=1$ and $t_{1}=0$, for any $p\geq 1$, we have
$$
E\left[\left|G^{\delta,i}(0, w(\ast))\right|^{2p}\right] \leq C_{p} \left\|E\left[\left|w(\ast)\right|^{2p}\right]\right\|_{C[0,\delta]}
= C_{p} \left\|t^{p}\right\|_{C[0,\delta]}\leq C_{p}\delta^{p}.
$$
\er

\br\label{r2} 
Since
$$
E\left[\left(\int_{k\delta}^{(k+1)\delta}\left|\dot{G}^{\delta,i}(s,w)\right| d s\right)^{p}\right]=E\left[\left(\int_{0}^{\delta}\left|\dot{G}^{\delta,i}(s,\theta_{k\delta} w)\right| d s\right)^{p}\right]=E\left[\left(\int_{0}^{\delta}\left|\dot{G}^{\delta,i}(s, w)\right| d s\right)^{p}\right]
$$
for $i=1,2,\cdots,r$, for any $0<\delta<1$, by H\"older's inequality, we have
\begin{equation}
E\left[\left(\int_{0}^{n_{1} \delta}\left|\dot{G}_{s}^{\delta,i_{1}}\right| d s\right)^{p_{1}}\cdots\left(\int_{0}^{n_{r} \delta}\left|\dot{G}_{s}^{\delta,i_{r}}\right| d s\right)^{p_{r}}\right]
\lesssim n_{1}^{p_{1}} n_{2}^{p_{2}} \cdots n_{r}^{p_{r}} \delta^{\frac{1}{2}\left(p_{1}+p_{2}+\cdots+p_{r}\right)}.
\end{equation}

\er

Let us introduce the following notations:
$$
s_{i j}(t, \delta) :=\frac{1}{2t} E\left[ \int_{0}^{t}G^{\delta,i}(s, w)\dot{G}^{\delta,j}(s, w)-G^{\delta,j}(s, w)\dot{G}^{\delta,i}(s, w)~ds\right],
$$
$$
c_{i j}^{*}(t, \delta) :=\frac{1}{\delta} E\left[\int_{0}^{\delta}\dot{G}^{\delta,i}(s, w)\left[G^{\delta,j}(t, w)-G^{\delta,j}(s, w)\right]ds\right]
$$
and
$$
c_{i j}(t, \delta) :=\frac{1}{t} E\left[\int_{0}^{t}\dot{G}^{\delta,i}(s, w)\left[G^{\delta,j}(t, w)-G^{\delta,j}(s, w)\right]ds\right].
$$

\br\label{r3}
We have
$
k c_{i j}(k\delta, \delta)= c_{i j}^{*}(k\delta, \delta)+(k-1)c_{i j}((k-1)\delta, \delta).
$
\er

\begin{assumption}\label{a1}
There exists a skew-symmetric $r \times r-$matrix $\left(s_{i j}\right)$ such that
$s_{i j}(\delta, \delta) \longrightarrow s_{i j}$ as $\delta\downarrow 0$.
Set
\begin{equation}\label{s_{i j}}
	c_{i j}=s_{i j}+\frac{1}{2} \delta_{i j}, \quad i, j=1,2, \ldots, r.
\end{equation}
\end{assumption}

We have the following result from $\cite{iw}$ (Section 7, Chapter VI).

\bp\label{p1}
Let $k(\delta) :(0,1] \rightarrow Z_{+}$ such that $k(\delta) \uparrow \infty$ as $\delta \downarrow 0$. Then
$$
\lim _{\delta \downarrow 0} c_{i j}(k(\delta) \delta, \delta)=c_{i j}.
$$
\ep

Now, we consider  stochastic differential equation with reflecting boundary $\partial D$,
\begin{equation}\label{i.1}
X_{t}=x+\int_{0}^{t}\sigma(X_{s})dW_{s}+\int_{0}^{t}\bar{b}(X_{s})ds+K_{t}
\end{equation}
and differential equation with reflecting boundary $\partial D$,
\begin{equation}\label{i,2}
X^{\delta}_{t}=x+\int_{0}^{t}\sigma(X^{\delta}_{s})dB^{\delta}_{s}+\int_{0}^{t} b(X^{\delta}_{s})ds+K^{\delta}_{t},
\end{equation}
where
$$
\bar{b}^{l}(X_{s})=b^{l}(X_{s})+\sum_{i,j=1}^{r}\sum_{\alpha=1}^{d}c_{ij}\sigma_{i}^{\alpha}\partial_{x_{\alpha}}\sigma_{j}^{l}(X_{s}),~for~l=1,\cdots,d,
$$
and $x\in\bar{D}$. The main result of this paper is the strong convergence
$$
E\left[\|X^{\delta,T}-X^{T}\|^p\right]\to 0~as~\delta\to 0,~\forall p>0.
$$

\section{Framework~and~some~results}

~~~~~The domain $D$, regarded as bounded now, is supposed to satisfy:

(A) There exists a constant $r_0>0$ such that for any $x\in\partial D$,
$$
\mathcal{N}_x=\mathcal{N}_{x,r_0}\neq\emptyset.
$$

(B) There exist constants $\delta>0$ and $\beta\geq1$ satisfying that for any $x\in\partial D$ there exists a unit vector $l_x$ such that
$$
\left\langle l_x,n\right\rangle \geq 1/\beta \quad \mbox{for any} \quad n\in\bigcup_{y\in B(x,\delta)\cap\partial D}\mathcal{N}_y,
$$ where $\left\langle \cdot,\cdot\right\rangle $ denotes the usual inner product in $\mathbb{R}^{d}$.

(C) There exists a $\mathcal{C}_b^2$ function $\varphi$ on $\mathbb{R}^d$ and a positive constant $\gamma$ such that for any $x\in\partial D$, $y\in\bar{D}$ and $n\in\mathcal{N}_x$,
$$
\left\langle y-x,n\right\rangle+\frac{1}{\gamma}\left\langle D\varphi (x), n\right\rangle|y-x|^2\geq0.
$$

(D) There exist $m \geq 1, \lambda>0, R>0, a_{1}, \ldots, a_{m} \in \mathbb{S}^{d-1}$ and $x_{1}, \ldots, x_{m} \in \partial D$
such that $\partial D \subset \bigcup_{i=1}^{m} B\left(x_{i}, R\right)$ and $x \in \partial D \cap B\left(x_{i}, 2 R\right) \Rightarrow n \cdot a_{i} \geq \lambda, \forall n \in \mathcal{N}_{x}$.
~\\

Under the above assumptions, the existence and uniqueness of $(\ref{i.1})$ and $(\ref{i,2})$ are proved by $\cite{ls}$. Here we have some results.

\bl\label{l1}
$($$\cite{ls}$, Theorem 1.1$)$. Assume (A),(D) hold, then for all $h\in C([0,\infty[,\mathbb{R}^{d})$ $(with ~h(0)\in\bar{D})$, there exists a unique solution $(x_{t},k_{t})$ of the $skorohod~problem$
$(h,D,\mathcal{N})$.
\el

\bl\label{l2} 
$($$\cite{aida}$, Lemma 2.4$)$. Assume condition (A) and the existence of the solution $(x,k)$ to the Skorohod problem
for a continuous bounded variation path $h$. Then the total variation of the solution $x$ has the
estimate:
$$
|x|_{s}^{t} \leq 2(\sqrt{2}+1)|h|_{s}^{t}.
$$
\el

\bl\label{l3} 
$($$\cite{stroock}$,  Theorem 3.5$)$. Assume $(D)$ holds and $(x,k)$ is a pair of solutions of the $skorohod~problem$
$(h,D,\mathcal{N})$, for any $0\leq s<t$, $\theta\in(0,\frac{1}{4})$,
$$
\left|k\right|_{0}^{t}-\left|k\right|_{0}^{s} \leq C\left((t-s) R^{-\frac{1}{\theta}}\left\|x\right\|_{[s, t],\theta}^{\frac{1}{\theta}}+1\right)\left\|k\right\|_{[s, t]},
$$
where $R$ is the constant given by condition (D).
\el

\bl\label{l4}
$($$\cite{aida}$, Lemma 2.3$)$. Assume $(A),(B)$ hold. Let $0<\theta \leq 1 .$ Then there exist positive constants $C_{1}, C_{2}, C_{3}$
which depend only on $\theta, \delta, \beta$ and $r$ in the Assumptions (A) and (B) such that
$$
|k|_{s}^{t}\leq C_{1}\left(1+\|h\|_{[s, t], \theta}^{C_{2}}(t-s)\right) e^{C_{3}\|h\|_{[s, t]}}\|h\|_{[s, t]}
$$
for all $0 \leq s<t \leq T$, where $(x(t),k(t))$ is a pair of solutions to $Skorohod~problem$ $(h,D,\mathcal{N})$.
\el

\section{Approximation ~theorem}

~~~~~We first state our main result.

\bt\label{t1} 
Under Assumption $\ref{a1}$, if $\sigma\in\mathcal{C}_{b}^{2}$ and $b$ is globally Lipschitz and bounded, then for any $p>0$, we have
$$
E\left[\|X^{\delta,T}-X^{T}\|^{p}\right]\to 0~as~\delta\to 0,
$$
where $X$, $\{X^{\delta}\}_{0<\delta<1}$ are solutions to $(\ref{i.1})$, $(\ref{i,2})$ respectively.
\et
The proof of Theorem $\ref{t1}$ is accomplished in $(\ref{theorem})$. Let $\mathcal{B}_{t}:=\sigma(W_{s}, s\leq t)$ and $\mathcal{F}_{t}:=\sigma({B}_{t}\cup\mathcal{N})$, where $\mathcal{N}$ denotes the $P$-negligible sets under $\mathcal{B}_{\infty}$.

\bl\label{l5} 
There exists a constant $C_{p}$ such that for any $t \in[0,T]$, $1\leq p$ and $0<\delta<1$.
$$
E\left[\left(|X^{\delta}|_{\underline{t}_{\delta}}^{t}\right)^{2 p}\right]\leq C_{p}\delta^{p},\quad E\left[\left(|K^{\delta}|_{\underline{t}_{\delta}}^{t}\right)^{2 p}\right]\leq C_{p}\delta^{p}.
$$
\el

\begin{proof}
Assume now
$$
L^{\delta}_{t}=x+\int_{0}^{t}\sigma(X^{\delta}_{s})d_{\delta}B_{s}+\int_{0}^{t} b(X^{\delta}_{s})ds.
$$
Since
$$
\left(|L^{\delta}|_{\underline{t}_{\delta}}^{t}\right)^{2 p}\lesssim \left(\int_{\underline{t}_{n}}^{\overline{t}_{\delta}}|\dot{B}^{\delta}_{s}|ds\right)^{2p}+(\overline{t}_{\delta}-\underline{t}_{\delta})^{2p},
$$
according to Definition (v) of $G^{\delta}_{t}$, we have
\begin{equation}
\begin{aligned}
E\left[\left(|L^{\delta}|_{\underline{t}_{\delta}}^{t}\right)^{2 p}\right]&\lesssim E\left[\left(\int_{\underline{t}_{\delta}}^{\overline{t}_{\delta}}|\dot{B}^{\delta}(s,w)|ds\right)^{2p}\right]+(\overline{t}_{\delta}-\underline{t}_{\delta})^{2p}\\
&= E\left[\left(\int_{0}^{2\delta}|\dot{G}^{\delta}(s,\theta_{\underline{\underline{t}}_{\delta}}w)|ds\right)^{2p}\right]+(\overline{t}_{\delta}-\underline{t}_{\delta})^{2p}\\
&=E\left[\left(\int_{0}^{2\delta}|\dot{G}^{\delta}(s,w)|ds\right)^{2p}\right]+(\overline{t}_{\delta}-\underline{t}_{\delta})^{2p}\\
&\lesssim~\delta^{p}.
\end{aligned}
\end{equation}
Using Lemma $\ref{l2}$, we get
$$
E\left[\left(|X^{\delta}|_{\underline{t}_{\delta}}^{t}\right)^{2 p}\right]\lesssim E\left[\left(|L^{\delta}|_{\underline{t}_{\delta}}^{t}\right)^{2 p}\right]\lesssim \delta^{p}
$$
and
$$
E\left[\left(|K^{\delta}|_{\underline{t}_{\delta}}^{t}\right)^{2 p}\right]\lesssim E\left[\left(|X^{\delta}|_{\underline{t}_{\delta}}^{t}\right)^{2 p}\right]+E\left[\left(|L^{\delta}|_{\underline{t}_{\delta}}^{t}\right)^{2 p}\right]\lesssim \delta^{p},
$$
which completes the prove.
\end{proof}

\bl\label{l6}
For any $s, t \in[0,T]$, $s<t$, we have
$$
E\left[\sup _{u, v \in[s, t]}\left|\int_{u}^{v} \sigma\left(X^{\delta}_{r}\right) dB^{\delta}_{r}\right|^{2p}\right]\leq C|t-s|^{p}
$$
and
$$
E\left[\left\|L^{\delta}\right\|_{[s, t]}^{2 p}\right]\leq C|t-s|^{p}.
$$
\el

\begin{proof}
If $\overline{s}_{\delta}<t_{\delta}$, write
\begin{equation}
	\int_{s}^{v} \sigma\left(X^{\delta}_{r}\right)dB^{\delta}_{r}
	=\int_{s}^{v} \sigma\left(X^{\delta}_{r}\right)-\sigma\left(X^{\delta}_{\underline{r}_{\delta}}\right)dB^{\delta}_{r}+\int_{s}^{v} \sigma\left(X^{\delta}_{\underline{r}_{\delta}}\right)dB^{\delta}_{r}
    :=\sum_{i=1}^{2}M_{i}(v)
\end{equation}
and
\begin{equation}
M_{2}(v)
=\int_{s}^{\overline{s}_{\delta}} \sigma\left(X^{\delta}_{\underline{r}_{\delta}}\right)dB^{\delta}_{r}
+\sum_{\overline{s}_{\delta}\leq k\delta\leq \underline{v}_{\delta}} \int_{k\delta}^{(k+1)\delta} \sigma\left(X^{\delta}_{\underline{r}_{\delta}}\right)dB^{\delta}_{r}
+\int_{v_{\delta}}^{v}\sigma\left(X^{\delta}_{\underline{r}_{\delta}}\right)dB^{\delta}_{r}.
\end{equation}

We have
\begin{equation}
\begin{aligned}
&\sum_{\overline{s}_{\delta}\leq k\delta\leq \underline{v}_{\delta}}\int_{k\delta}^{(k+1)\delta} \sigma\left(X^{\delta}_{\underline{r}_{\delta}}\right)dB^{\delta}_{r}\\
=&\sum_{\overline{s}_{\delta}\leq k\delta\leq \underline{v}_{\delta}} \sigma\left(X^{\delta}(k\delta-\delta,w)\right)\big[B^{\delta}(k\delta+\delta,w)-B^{\delta}(k\delta,w)\big]\\
=&\sum_{\overline{s}_{\delta}\leq k\delta\leq \underline{v}_{\delta}} \sigma\left(X^{\delta}(k\delta-\delta,w)\right)\big[W(k\delta)-W(k\delta-\delta)\big]\\
&+\sum_{\overline{s}_{\delta}\leq k\delta\leq \underline{v}_{\delta}} \sigma\left(X^{\delta}(k\delta-\delta,w)\right)\left(G^{\delta}(0,\theta_{k\delta}w)-G^{\delta}(0,\theta_{(k-1)\delta}w)\right)\\
:=&\sum_{i=1}^{2}H_{i}(v).
\end{aligned}
\end{equation}
Obviously, it holds that
$$
\sum_{\overline{s}_{\delta}\leq k\delta\leq \underline{v}_{\delta}} \sigma\left(X^{\delta}(k\delta-\delta,w)\right)\big[W(k\delta)-W(k\delta-\delta)\big]
=\int_{s_{\delta}}^{\underline{v}_{\delta}}\sigma\left(X^{\delta}_{r_{\delta}}\right) dW_{r},
$$
thus BDG's inequality implies
$$
E\left[\sup _{v \in[s, t]}\left|H_{1}(v)\right|^{2p }\right]\leq C|\underline{t}_{\delta}-s_{\delta}|^{p}.
$$
Note that
$$
S'_{n}=\sum_{k=1}^{n} \sigma\left(X^{\delta}(k\delta-\delta,w)\right)G^{\delta}(0,\theta_{k\delta}w)
$$
is a $\{\mathcal{F}_{n}\}-$martingale where $\mathcal{F}_{n}=\mathcal{B}_{(n+1)\delta}$, and
$$
S''_{n}=-\sum_{k=1}^{n} \sigma\left(X^{\delta}(k\delta-\delta,w)\right)G^{\delta}(0,\theta_{(k-1)\delta}w)
$$
is a $\{\mathcal{H}_{n}\}-$martingale where $\mathcal{F}_{n}=\mathcal{B}_{n\delta}$. Since
$$
E\left[\max_{\overline{s}_{\delta}\leq k\delta\leq \underline{v}_{\delta}}|\sigma\left(X^{\delta}(k\delta-\delta,w)\right)G^{\delta}(0,\theta_{k\delta}w)|^{2p}\right]
\lesssim \frac{|t_{\delta}-\overline{s}_{\delta}|}{\delta}\delta^{p}
\leq |t-s|^{p}
$$
and
$$
E\left[|\sigma\left(X^{\delta}(s_{\delta},w)\right)G^{\delta}(0,\theta_{s_{\delta}}w)|^{2p}\right]\lesssim \delta^{p}\leq |t-s|^{p},
$$
by martingale inequality, it's easy to see that
\begin{equation}
\begin{aligned}
&E\left[\sup _{v \in[s, t]}\left|H_{2}(v)\right|^{2 p}\right]\\
\lesssim &E\left[~\left|\sum_{k=\overline{s}_{\delta}/\delta}^{t_{\delta}/\delta} \sigma\left(X^{\delta}(k\delta-\delta,w)\right)\left(G^{\delta}(0,\theta_{k\delta}w)-G^{\delta}(0,\theta_{(k-1)\delta}w)\right)\right|^{2p }\right]+|t-s|^{p}\\
=& E\left[~\left|G^{\delta}\left(0,\sum_{k=\overline{s}_{\delta}/\delta}^{t_{\delta}/\delta} \sigma\left(X^{\delta}(k\delta-\delta,w)\right)(\theta_{k\delta}w-\theta_{(k-1)\delta}w)\right)\right|^{2p }~\right]+|t-s|^{p}\\
\lesssim & \left\|E\left[~\left|\sum_{k=\overline{s}_{\delta}/\delta}^{t_{\delta}/\delta} \sigma(X^{\delta}(k\delta-\delta,w))(w(k\delta+\ast)-w((k-1)\delta+\ast))\right|^{2p}~\right]\right\|_{C[0,\delta]}\\
&+
\left\|E\left[~\left|\sum_{k=\overline{s}_{\delta}/\delta}^{t_{\delta}/\delta} \sigma(X^{\delta}(k\delta-\delta,w))(w(k\delta)-w((k-1)\delta)\right|^{2p}~\right]\right\|_{C[0,\delta]}+|t-s|^{p}\\
= & \left\|E\left[~\left|\int_{s_{\delta}}^{t_{\delta}}\sigma(X^{\delta}(r_{\delta},w))dw_{r+\ast}\right|^{2p}~\right]\right\|_{C[0,\delta]}+\left\|E\left[~\left|\int_{s_{\delta}}^{t_{\delta}}\sigma(X^{\delta}(r_{\delta},w))dw_{r}\right|^{2p}~\right]\right\|_{C[0,\delta]}+|t-s|^{p}\\
\lesssim & \left\|(t_{\delta}-s_{\delta})^{p}\right\|_{C[0,\delta]}+|t-s|^{p}\leq|t-s|^{p}.
\end{aligned}
\end{equation}
Moreover
\begin{equation}
\begin{aligned}
E\left[\sup _{v\in[s, t]}\left|\int_{v_{\delta}}^{v} \sigma\left(X^{\delta}_{\underline{r}_{\delta}}\right)dB^{\delta}_{r}\right|^{2p}\right]
\leq &~E\left[\sup _{v\in[s, t]}\left(\int_{v_{\delta}}^{\overline{v}_{\delta}} |\dot{B}^{\delta}_{r}|dr \right)^{2p}\right]
\leq \left(E\left[\sup _{v\in[s, t]}\left(\int_{v_{n}}^{\overline{v}_{\delta}} |\dot{B}^{\delta}_{r}|dr \right)^{4p}\right]\right)^{\frac{1}{2}}\\
\leq&
\left(\sum_{s_{\delta}\leq k\delta\leq t_{\delta}}E\left[\left(\int_{k\delta}^{(k+1)\delta} |\dot{B}^{\delta}_{r}|dr \right)^{4p}\right]\right)^{\frac{1}{2}}
\lesssim\left(\frac{|t_{\delta}-\underline{s}_{\delta}|}{\delta}\delta^{2p}\right)^{\frac{1}{2}}\\
\lesssim &~~|t-s|^{p},
\end{aligned}
\end{equation}
and we can also get
$$
E\left[\left|\int_{s}^{\overline{s}_{\delta}}\sigma\left(X^{\delta}_{\underline{r}_{\delta}}\right)dB^{\delta}_{r}\right|^{2p}\right]\lesssim|t-s|^{p}
$$
in the same way. All the estimates before show
$$
E\left[\sup _{v \in[s, t]}\left|M_{2}(v)\right|^{2p}\right]\lesssim|t-s|^{p}.
$$
As for the term $M_{1}$, we have
\begin{equation}
\begin{aligned}
E\left[\sup _{v \in[s, t]}\left|M_{1}(v)\right|^{2p}\right]
\lesssim & E\left[\left(\int_{\underline{s}_{\delta}}^{\overline{t}_{\delta}}|X^{\delta}|_{\underline{r}_{\delta}}^{\overline{r}_{\delta}}|\dot{B}^{\delta}_{r}|dr\right)^{2p}\right]
\\
\leq &\left(\frac{|t_{\delta}-\underline{s}_{\delta}|}{\delta}\right)^{2p-1}\sum_{k=\underline{s}_{\delta}/\delta}^{t_{\delta}/\delta}E\left[\left(|X^{\delta}|_{(k-1)\delta}^{(k+1)\delta}\int_{k\delta}^{(k+1)\delta}|\dot{B}^{\delta}_{r}|dr\right)^{2p}\right]\\
\lesssim & \left(\frac{|t_{\delta}-\underline{s}_{\delta}|}{\delta}\right)^{2p}\delta^{2p}
\lesssim~|t-s|^{2p}.
\end{aligned}
\end{equation}
Then
\begin{equation}
\begin{aligned}
	E\left[\sup _{u, v \in[s, t]}\left|\int_{u}^{v} \sigma\left(X^{\delta}_{r}\right)dB^{\delta}_{r}\right|^{2p}\right]
	&\lesssim\sum_{i=1}^{2}E\left[\sup _{v \in[s, t]}\left|M_{i}(v)\right|^{2p}\right]
	\lesssim |t-s|^{p}.
\end{aligned}
\end{equation}
If $\overline{s}_{\delta}\geq t_{\delta}$,  the result is trivial by Definition (v). The second inequality is an immediate result from the first one, thus we complete the prove.
\end{proof}

\bl\label{l7} 
For any $s,t\in[0,T]$, $s<t$ and $p\geq1$, there is a constant $C_{p}$ such that
$$
E\left[\|X^{\delta}\|_{[s,t]}^{4p}\right]\leq C_{p}|t-s|^{p},~E\left[\|K^{\delta}\|_{[s,t]}^{4p}\right]\leq C_{p}|t-s|^{p}.
$$
\el

\begin{proof}
If $\overline{s}_{\delta}\geq t_{\delta}$, a similar argument in Lemma $\ref{l5}$ leads to the result. Otherwise, for all $r\in[s,t]$, we have
\begin{equation}
\begin{aligned}
&e^{-\frac{2}{\gamma}\varphi(X^{\delta}_{r})}|X^{\delta}_{r}-X^{\delta}_{s}|^{2}\\
=&2\int_{s}^{r}e^{-\frac{2}{\gamma}\varphi(X^{\delta}_{u})}(X^{\delta}_{u}-X^{\delta}_{s})^{*}\sigma(X^{\delta}_{u})dB^{\delta}_{u}\\
+&2\int_{s}^{r}e^{-\frac{2}{\gamma}\varphi( X^{\delta}_{u})}(X^{\delta}_{u}-X^{\delta}_{s})^{*}b(X^{\delta}_{u})du\\
+&2\int_{s}^{r}e^{-\frac{2}{\gamma}\varphi(X^{\delta}_{u})}( X^{\delta}_{u}-X^{\delta}_{s})^{*}dK^{\delta}_{u}\\
-&\frac{2}{\gamma}\int_{s}^{r}e^{-\frac{2}{\gamma}\varphi(X^{\delta}_{u})}|X^{\delta}_{u}-X^{\delta}_{s}|^{2}D\varphi(X^{\delta}_{u})\sigma(X^{\delta}_{u})dB^{\delta}_{u}\\
-&\frac{2}{\gamma}\int_{s}^{r}e^{-\frac{2}{\gamma}\varphi(X^{\delta}_{u})}|X^{\delta}_{u}-X^{\delta}_{s}|^{2}D\varphi(X^{\delta}_{u})b(X^{\delta}_{u})du\\
-&\frac{2}{\gamma}\int_{s}^{r}e^{-\frac{2}{\gamma}\varphi(X^{\delta}_{u})}|X^{\delta}_{u}-X^{\delta}_{s}|^{2}D\varphi(X^{\delta}_{u})dK^{\delta}_{u}\\
:=&\sum_{i=1}^{6}U_{i}(r).
\end{aligned}
\end{equation}
It's obvious that
$$
U_{3}(r)+U_{6}(r)\leq 0.
$$
Since $D$ is bounded, it's easy to see that
$$
E\left[\sup_{s\leq v\leq r}|U_{2}(v)|^{2p}\right]\lesssim |r-s|^{2p}
$$
and
$$
E\left[\sup_{s\leq v\leq r}|U_{5}(v)|^{2p}\right]\lesssim |r-s|^{2p}.
$$
As for $U_{1}(r)$, write
\begin{equation}
\begin{aligned}
&~U_{1}(r)\\
=&2\int_{s}^{r}e^{-\frac{2}{\gamma}\varphi(X^{\delta}_{\underline{u}_{\delta}})}(X^{\delta}_{\underline{u}_{\delta}}-X^{\delta}_{s})^{*}\sigma(X^{\delta}_{\underline{u}_{\delta}})dB^{\delta}_{u}\\
+&2\int_{s}^{r}e^{-\frac{2}{\gamma}\varphi(X^{\delta}_{u})}X^{\delta,*}_{u}\sigma(X^{\delta}_{u})-e^{-\frac{2}{\gamma}\varphi(X^{\delta}_{\underline{u}_{\delta}})}X^{\delta,*}_{\underline{u}_{\delta}}\sigma(X^{\delta}_{\underline{u}_{\delta}})dB^{\delta}_{u}\\
-&2\int_{s}^{r}X^{\delta,*}_{s}\left(e^{-\frac{2}{\gamma}\varphi(X^{\delta}_{u})}\sigma(X^{\delta}_{u})-e^{-\frac{2}{\gamma}\varphi(X^{\delta}_{\underline{u}_{\delta}})}\sigma(X^{\delta}_{\underline{u}_{\delta}})\right)dB^{\delta}_{u}\\
:=&\sum_{i=1}^{3}U_{1i}(r).
\end{aligned}
\end{equation}
Obviously, for $i=2,3$, we have
\begin{equation}
\begin{aligned}
&E\left[\sup_{s\leq v\leq r}|U_{1i}(v)|^{2p}\right]\\
\lesssim & (r-s)^{2p-1}E\left[\sum_{k=s_{\delta}/\delta}^{r_{\delta}/\delta}\int_{k\delta}^{(k+1)\delta}\left(|X^{\delta}|_{(k-1)\delta}^{(k+1)\delta}|\dot{B}^{\delta}_{u}|\right)^{2p}du\right]\\
\leq & (r-s)^{2p-1}\sum_{k=s_{\delta}/\delta}^{r_{\delta}/\delta}\left(E\left[\left(|X^{\delta}|_{(k-1)\delta}^{(k+1)\delta}\right)^{4p}\right]\right)^{\frac{1}{2}}\left(E\left[\left(\int_{k\delta}^{(k+1)\delta}|\dot{B}^{\delta}_{u}|^{2p}du\right)^{2}\right]\right)^{\frac{1}{2}}\\
\leq & (r-s)^{2p-1}\frac{\overline{r}_{\delta}-s_{\delta}}{\delta}\delta\lesssim (r-s)^{2p-1}.
\end{aligned}
\end{equation}

Write
\begin{equation}
\begin{aligned}
&\int_{\overline{s}_{\delta}}^{r_{\delta}}e^{-\frac{2}{\gamma}\varphi(X^{\delta}_{\underline{u}_{\delta}})}(X^{\delta}_{\underline{u}_{\delta}}-X^{\delta}_{s})^{*}\sigma(X^{\delta}_{\underline{u}_{\delta}})dB^{\delta}_{u}\\
=&\sum_{k=\overline{s}_{\delta}}^{\underline{r}_{\delta}}e^{-\frac{2}{\gamma}\varphi(X^{\delta}_{(k-1)\delta})}(X^{\delta}_{(k-1)\delta}-X^{\delta}_{s})^{*}\sigma(X^{\delta}_{(k-1)\delta})\left(w(k\delta)-w(k\delta-\delta)\right)\\
&+\sum_{k=\overline{s}_{\delta}}^{\underline{r}_{\delta}}e^{-\frac{2}{\gamma}\varphi(X^{\delta}_{(k-1)\delta})}(X^{\delta}_{(k-1)\delta}-X^{\delta}_{s})^{*}\sigma(X^{\delta}_{(k-1)\delta})\left(G^{\delta}(0,\theta_{k\delta}w)-G^{\delta}(0,\theta_{(k-1)\delta}w)\right)\\
:=&\sum_{i=1}^{2}U_{11i}(r).
\end{aligned}
\end{equation}
Note that $D$ is bounded, BDG's inequality implies
$$
E\left[\sup_{s\leq v\leq r}|U_{111}(v)|^{2p}\right]\lesssim (r-s)^{p-1}\int_{s_{\delta}}^{r_{\delta}}E\left[|X^{\delta}_{u_{\delta}}-X^{\delta}_{s}|^{2p}\right]du
\lesssim (r-s)^{p},
$$
and a similar argument as we handle $H_{2}$ in Lemma $\ref{l6}$ leads to
$$
E\left[\sup_{s\leq v\leq r}|U_{112}(v)|^{2p}\right]\lesssim (r-s)^{p}.
$$
Thus we claim
$$
E\left[\sup_{s\leq v\leq r}|U_{1}(v)|^{2p}\right]\lesssim (r-s)^{p}.
$$
Now we consider $U_{4}(r)$,
\begin{equation}
\begin{aligned}
-\frac{2}{\gamma}&\int_{s}^{r}e^{-\frac{2}{\gamma}\varphi(X^{\delta}_{u})}|X^{\delta}_{u}-X^{\delta}_{s}|^{2}D\varphi(X^{\delta}_{u})\sigma(X^{\delta}_{u})dB^{\delta}_{u}\\
=-\frac{2}{\gamma}&\int_{s}^{r}e^{-\frac{2}{\gamma}\varphi(X^{\delta}_{\underline{u}_{\delta}})}|X^{\delta}_{u}-X^{\delta}_{\underline{u}_{\delta}}|^{2}D\varphi(X^{\delta}_{\underline{u}_{\delta}})\sigma(X^{\delta}_{\underline{u}_{\delta}})dB^{\delta}_{u}\\
-\frac{4}{\gamma}&\int_{s}^{r}e^{-\frac{2}{\gamma}\varphi(X^{\delta}_{\underline{u}_{\delta}})}(X^{\delta}_{u}-X^{\delta}_{\underline{u}_{\delta}})(X^{\delta}_{\underline{u}_{\delta}}-X^{\delta}_{s})D\varphi(X^{\delta}_{\underline{u}_{\delta}})\sigma(X^{\delta}_{\underline{u}_{\delta}})dB^{\delta}_{u}\\
-\frac{2}{\gamma}&\int_{s}^{r}e^{-\frac{2}{\gamma}\varphi(X^{\delta}_{\underline{u}_{\delta}})}|X^{\delta}_{\underline{u}_{\delta}}-X^{\delta}_{s}|^{2}D\varphi(X^{\delta}_{\underline{u}_{\delta}})\sigma(X^{\delta}_{\underline{u}_{\delta}})dB^{\delta}_{u}\\
-\frac{2}{\gamma}&\int_{s}^{r}|X^{\delta}_{u}-X^{\delta}_{s}|^{2}\left(e^{-\frac{2}{\gamma}\varphi(X^{\delta}_{u})}D\varphi(X^{\delta}_{u})\sigma(X^{\delta}_{u})-e^{-\frac{2}{\gamma}\varphi(X^{\delta}_{\underline{u}_{\delta}})}D\varphi(X^{\delta}_{\underline{u}_{\delta}})\sigma(X^{\delta}_{\underline{u}_{\delta}})\right)dB^{\delta}_{u}\\
:=\sum_{i=1}^{4}&U_{4i}(r).
\end{aligned}
\end{equation}
Since
$$
E\left[\left(|X^{\delta}|_{\underline{t}_{\delta}}^{t}\right)^{2 p}\right]\leq C_{p}\delta^{p},
$$
reasoning in the same way as we do to $U_{12}$, we have
$$
E\left[\sup_{s\leq v\leq r}|U_{4i}(v)|^{2p}\right]\lesssim |r-s|^{p},
$$
for $i=1,2,4$. Also
\begin{equation}
\begin{aligned}
-\frac{2}{\gamma}&\int_{\overline{s}_{\delta}}^{r_{\delta}}e^{-\frac{2}{\gamma}\varphi(X^{\delta}_{\underline{u}_{\delta}})}|X^{\delta}_{\underline{u}_{\delta}}-X^{\delta}_{s}|^{2}D\varphi(X^{\delta}_{\underline{u}_{\delta}})\sigma(X^{\delta}_{\underline{u}_{\delta}})dB^{\delta}_{u}\\
=-\frac{2}{\gamma}&\sum_{k=\overline{s}_{\delta}}^{\underline{r}_{\delta}}e^{-\frac{2}{\gamma}\varphi(X^{\delta}_{(k-1)\delta})}|X^{\delta}_{(k-1)\delta}-X^{\delta}_{s}|^{2}D\varphi(X^{\delta}_{(k-1)\delta})\sigma(X^{\delta}_{(k-1)\delta})\left(w(k\delta)-w(k\delta-\delta)\right)\\
-\frac{2}{\gamma}&\sum_{k=\overline{s}_{\delta}}^{\underline{r}_{\delta}}e^{-\frac{2}{\gamma}\varphi(X^{\delta}_{(k-1)\delta})}|X^{\delta}_{(k-1)\delta}-X^{\delta}_{s}|^{2}D\varphi(X^{\delta}_{(k-1)\delta})\sigma(X^{\delta}_{(k-1)\delta})\left(_{\delta}G(0,\theta_{k\delta}w)-_{\delta}G(0,\theta_{(k-1)\delta}w)\right).
\end{aligned}
\end{equation}
A same way as we do to $U_{11}$ shows
$$
E\left[\sup_{s\leq v\leq r}|U_{43}(v)|^{2p}\right]\lesssim (r-s)^{p}.
$$
Summing up all above, we claim
$$
E\left[\sup_{s\leq r\leq t}|X^{\delta}_{r}-X^{\delta}_{s}|^{4p}\right]\lesssim (r-s)^{p}.
$$
Consequently
$$
E\left[\|X^{\delta}\|_{[s,t]}^{4p}\right]\lesssim (r-s)^{p}.
$$
Also, by Lemma $\ref{l6}$, we get
$$
E\left[\|K^{\delta}\|_{[s,t]}^{4p}\right]\lesssim E\left[\|X^{\delta}\|_{[s,t]}^{4p}\right]+E\left[\left\|L^{\delta}\right\|_{[s, t]}^{4p}\right]\lesssim(r-s)^{p}.
$$
\end{proof}

\br\label{r4}
By Lemma $\ref{l7}$ and $Kolmogorov's ~continuity~criterion$, $p\geq 1$, for any $\theta\in\left(0, \frac{1}{4}\right)$, we have
$$
E\left[\|X^{\delta}\|_{[s,t],\theta}^{p}\right]\leq C_{p,\theta},~E\left[\|K^{\delta}\|_{[s,t],\theta}^{p}\right]\leq C_{p,\theta}.
$$
\er

\br\label{r5} 
By Remark $\ref{r4}$ and Lemma $\ref{l3}$, $p\geq 1$, we have
$$
E\left[\left(\left|K^{\delta}\right|_{0}^{T}\right)^{p}\right]<C_{p}.
$$
Since $\forall\theta\in(0,\frac{1}{4})$, there is a $p$ large enough satisfying $\theta\in(0,\frac{1}{4}-\frac{1}{4p})$,
$$
E\left[\sup_{0\leq u,v\leq T}|K^{\delta}_{u}-K^{\delta}_{v}|^{4p}\right]\lesssim E\left[\sup_{0\leq u,v\leq T,u\neq v}\left(\frac{|K^{\delta}_{u}-K^{\delta}_{v}|}{|u-v|^{\theta}}\right)^{4p}\right]<\infty.
$$
\er
~\\

Following notations in $\cite{iw}$, let $n(\delta) :(0,1] \longrightarrow Z_{+}$ such that $n(\delta) \uparrow \infty~as~\delta \downarrow 0$ and $n(\delta)^{5} \delta \downarrow 0~as~\delta \downarrow 0$. Assume that
\begin{equation}
\begin{aligned}
\widetilde{\delta}&=n(\delta) \delta,\\
[s](\widetilde{\delta})&=k\widetilde{\delta},\\
[s]^{+}(\widetilde{\delta})&=(k+1) \widetilde{\delta},\\
[s]^{-}(\widetilde{\delta})&=\left((k-1)\vee 0\right) \widetilde{\delta},
\end{aligned}
\end{equation}
if $k\widetilde{\delta}\leq s<(k+1)\widetilde{\delta}$, and
$$
S_{\delta}(t)=[t](\widetilde{\delta})/\widetilde{\delta}.
$$

\bl\label{l8} 
Let
$$
\mu_{\delta}(t)=e^{-\frac{2}{\gamma}\left(\varphi\left(X_{t}\right)+\varphi\left(X^{\delta}_{t}\right)\right)},~ m_{\delta}(t)=e^{-\frac{2}{\gamma}\left(\varphi\left(X_{t}\right)+\varphi\left(X^{\delta}_{t}\right)\right)}\left|X_{t}-X^{\delta}_{t}\right|^{2}.
$$
Then
$$
E\left[|\mu_{\delta}(t)-\mu_{\delta}([t](\widetilde{\delta}))|^{2}\right]\lesssim\widetilde{\delta}^{\frac{1}{2}},~E\left[|m_{\delta}(t)-m_{\delta}([t](\widetilde{\delta}))|^{2}\right]\lesssim\widetilde{\delta}^{\frac{1}{2}},~\forall t\in[0,T],
$$
where $\varphi$ is given by condition $(C)$ and $X$, $\{X^{\delta}\}_{0<\delta<1}$ is a sequence of solutions to $(\ref{i.1})$ and $(\ref{i,2})$.
\el

\begin{proof}
Since $\varphi\in C_{b}^{2}$, by Lemma $\ref{l4}$ and Lemma $\ref{l7}$, we have
$$
E\left[|\mu_{\delta}(t)-\mu_{\delta}([t](\widetilde{\delta}))|^{2}\right]
\lesssim E\left[|X_{t}-X_{[t](\widetilde{\delta})}|^{2}\right]+E\left[|X^{\delta}_{t}-X^{\delta}_{[t](\widetilde{\delta})}|^{2}\right]
\lesssim\widetilde{\delta}^{\frac{1}{2}}.
$$
By virtue of the boundness of $D$, we get
\begin{equation}
\begin{aligned}
&E\left[|m_{\delta}(t)-m_{\delta}([t](\widetilde{\delta}))|^{2}\right]\\
\lesssim~&E\left[|\mu_{\delta}(t)-\mu_{\delta}([t](\widetilde{\delta}))|^{2}\right]+E\left[|X_{t}-X_{[t](\widetilde{\delta})}|^{2}\right]+E\left[|X^{\delta}_{t}-X^{\delta}_{[t](\widetilde{\delta})}|^{2}\right]\\
\lesssim~&\widetilde{\delta}^{\frac{1}{2}}.
\end{aligned}
\end{equation}
which is the desired result.
\end{proof}

\bp\label{p2}
Under Assumption $\ref{a1}$,  if $\sigma\in\mathcal{C}_{b}^{2}$ and $b$ is globally Lipschitz and bounded, then
$$
E\left[\sup_{u\leq t}|X_{u}-X^{\delta}_{u}|^4\right]\lesssim\int_{0}^{t}E[\sup_{u\leq s}|X_{u}-X^{\delta}_{u}|^{4}]ds+o(1),~\forall t\in[0,T],
$$
 where
$$
o(1)=n(\delta)^{\frac{5}{2}}\delta^{\frac{1}{2}}+n(\delta)^{-1}+|(c_{ij})(\widetilde{\delta},\delta)-(c_{ij})|^{2}
$$
and
$X$, $\{X^{\delta}\}_{0<\delta<1}$ are solutions to $(\ref{i.1})$, $(\ref{i,2})$ respectively. Convention: We say a process $F_{t}$ is trivial if $E[\sup_{t\leq T}|F_{t}|^{2}]\lesssim o(1)$.
\ep

\begin{proof}
By It\^o's formula, we have
\begin{equation}
\begin{aligned}
	&e^{-\frac{2}{\gamma}(\varphi(X_{t})+\varphi(X^{\delta}_{t}))}|X_{t}-X^{\delta}_{t}|^{2}\\
	=&-\frac{2}{\gamma}\int_{0}^{t}m_{\delta}(u)I_{\delta}(u)du\\
	&-\frac{2}{\gamma}\int_{0}^{t}m_{\delta}(u)D\varphi(X^{\delta}_{u})\sigma(X^{\delta}_{u})dB^{\delta}_{u}\\
	&-\frac{2}{\gamma}\int_{0}^{t}m_{\delta}(u)D\varphi(X_{u})\sigma(X_{u})dW_{u}
\end{aligned}
\end{equation}	
\begin{equation}
\begin{aligned}	
	&+\int_{0}^{t}\mu_{\delta}(u)[-\frac{2}{\gamma}|X_{u}-X^{\delta}_{u}|^{2}( D\varphi(X^{\delta}_{u})dK^{\delta}_{u}+ D\varphi(X_{u})dK_{u})+2( X_{u}-X^{\delta}_{u})^{*}(dK_{u}-dK^{\delta}_{u})]\\
	&+2\int_{0}^{t}\mu_{\delta}(u)(X_{u}-X^{\delta}_{u})^{*}(b(X_{u})-b(X^{\delta}_{u}))du\\
	&-\frac{4}{\gamma}\int_{0}^{t}\mu_{\delta}(u)D\varphi(X_{u})\sigma(X_{u})\sigma^{*}(X_{u})(X_{u}-X^{\delta}_{u})du\\
	&-2\int_{0}^{t}\mu_{\delta}(u)(X_{u}-X^{\delta}_{u})^{*}\sigma(X^{\delta}_{u})dB^{\delta}_{u}\\
	&+2\int_{0}^{t}\mu_{\delta}(u)(X_{u}-X^{\delta}_{u})^{*}\sigma(X_{u})dW_u\\
	&+2\int_{0}^{t}\mu_{\delta}(u)(X_{u}-X^{\delta}_{u})^{*}\left(\bar{b}-b\right)(X_{u})du\\
	&+\int_{0}^{t}\mu_{\delta}(u)tr(\sigma^{*}\sigma(X_{u}))du\\
	:=&\sum_{i=1}^{10}I_{i}(t),
\end{aligned}
\end{equation}
and
\begin{equation}
	I_{\delta}(u)=D\varphi(X^{\delta}_{u})b(X^{\delta}_{u})+D\varphi(X_{u})\bar{b}(X_{u})
	+\frac{1}{2}tr(\sigma^{*}D^{2}\varphi\sigma(X_{u}))-\frac{1}{\gamma}|D\varphi(X_{u})b(X_{u})|^{2}.
\end{equation}
It's easy to see that
$$
E[\sup_{s\leq t}|I_{1}(s)|^{2}]\lesssim\int_{0}^{t}E[m_{\delta}^{2}(u)]du.
$$
BDG's inequality implies
$$
E[\sup_{s\leq t}|I_{3}(s)|^{2}]\lesssim\int_{0}^{t}E[m_{\delta}^{2}(u)]du,
$$
and Condition (C) yields
$$
I_{4}(t)\leq 0.
$$
Since $\varphi\in \mathcal{C}_{b}^{2}$, we get
$$
E[\sup_{s\leq t}|I_{5}(s)|^{2}]\lesssim\int_{0}^{t}E[m_{\delta}^{2}(u)]du.
$$
As for the term $I_{2}$,
\begin{equation}
\begin{aligned}
&~~I_{2}(t)\\
=&-\frac{2}{\gamma}\int_{[t](\widetilde{\delta})}^{t}m_{\delta}(u)D\varphi(X^{\delta}_{u})\sigma(X^{\delta}_{u})dB^{\delta}_{u}\\
&-\frac{2}{\gamma}\int_{\widetilde{\delta}}^{[t](\widetilde{\delta})}m_{\delta}(u)D\varphi(X^{\delta}_{u})\sigma(X^{\delta}_{u})dB^{\delta}_{u}\\
&-\frac{2}{\gamma}\int_{0}^{\widetilde{\delta}}m_{\delta}(u)D\varphi(X^{\delta}_{u})\sigma(X^{\delta}_{u})dB^{\delta}_{u}\\
:=&\sum_{i=1}^{3}I_{2i}(t).
\end{aligned}
\end{equation}
It follows that
\begin{equation}
\begin{aligned}
E\left[\sup_{s\leq t}|I_{21}(s)|^{2}\right]\leq&\left(E\left[\sup_{s\in[0,T]}|I_{21}(s)|^{4}\right]\right)^{\frac{1}{2}}\lesssim\left(\sum_{k}E\left[\left(\int_{k\widetilde{\delta}}^{(k+1)\widetilde{\delta}}|\dot{B}^{\delta}_{u}|du\right)^{4}\right]\right)^{\frac{1}{2}}\\
\lesssim~&\left(\frac{1}{n(\delta)\delta}n(\delta)^{4}\delta^{2}\right)^{\frac{1}{2}}
\leq~n(\delta)^{\frac{3}{2}}\delta^{\frac{1}{2}},
\end{aligned}
\end{equation}
and a similar argument leads to
$$
E\left[\sup_{s\leq t}|I_{23}(s)|^{2}\right]\lesssim n(\delta)^{\frac{3}{2}}\delta^{\frac{1}{2}}.
$$

Now split $I_{22}(t)$ as
\begin{equation}
\begin{aligned}
&I_{22}(t)\\
=&-\frac{2}{\gamma}\sum_{k=1}^{S_{\delta}(t)-1}\int_{k\widetilde{\delta}}^{(k+1)\widetilde{\delta}}m_{\delta}(u)D\varphi\sigma(X^{\delta}_{u})dB^{\delta}_{u}\\
=&-\frac{2}{\gamma}\sum_{k=1}^{S_{\delta}(t)-1}\int_{k\widetilde{\delta}}^{(k+1)\widetilde{\delta}}m_{\delta}(u)\left(D\varphi\sigma(X^{\delta}_{u})-D\varphi\sigma(X^{\delta}_{k\widetilde{\delta}-\delta})\right)dB^{\delta}_{u}\\
&-\frac{2}{\gamma}\sum_{k=1}^{S_{\delta}(t)-1}\int_{k\widetilde{\delta}}^{(k+1)\widetilde{\delta}}\left(m_{\delta}(u)-m_{\delta}(k\widetilde{\delta}-\delta)\right)D\varphi\sigma(X^{\delta}_{k\widetilde{\delta}-\delta})dB^{\delta}_{u}\\
&-\frac{2}{\gamma}\sum_{k=1}^{S_{\delta}(t)-1}\int_{k\widetilde{\delta}}^{(k+1)\widetilde{\delta}}m_{\delta}(k\widetilde{\delta}-\delta)D\varphi\sigma(X^{\delta}_{k\widetilde{\delta}-\delta})dB^{\delta}_{u}\\
:=&\sum_{i=1}^{3}\sum_{k=1}^{S_{\delta}(t)-1}I_{22i}^{k}.
\end{aligned}
\end{equation}
We have
\begin{equation}
\begin{aligned}
&~~~~I_{221}^{k}\\
=&-\frac{2}{\gamma}\sum_{i,j=1}^{r}\int_{k\widetilde{\delta}}^{(k+1)\widetilde{\delta}}m_{\delta}(k\widetilde{\delta}-\delta)\int_{k\widetilde{\delta}-\delta}^{u}(D(D\varphi\sigma)_{i}\sigma)_{j}(X^{\delta}_{k\widetilde{\delta}-\delta})dB^{\delta,j}_{s}dB^{\delta,i}_{u}\\
&-\frac{2}{\gamma}\sum_{i,j=1}^{r}\int_{k\widetilde{\delta}}^{(k+1)\widetilde{\delta}}m_{\delta}(u)\int_{k\widetilde{\delta}-\delta}^{u}(D(D\varphi\sigma)_{i}\sigma)_{j}(X^{\delta}_{u})-(D(D\varphi\sigma)_{i}\sigma)_{j}(X^{\delta}_{k\widetilde{\delta}-\delta})dB^{\delta,j}_{s}dB^{\delta,i}_{u}\\
&-\frac{2}{\gamma}\sum_{i,j=1}^{r}\int_{k\widetilde{\delta}}^{(k+1)\widetilde{\delta}}(m_{\delta}(u)-m_{\delta}(k\widetilde{\delta}-\delta))\int_{k\widetilde{\delta}-\delta}^{u}(D(D\varphi\sigma)_{i}\sigma)_{j}(X^{\delta}_{k\widetilde{\delta}-\delta})dB^{\delta,j}_{s}dB^{\delta,i}_{u}\\
&-\frac{2}{\gamma}\sum_{i=1}^{r}\int_{k\widetilde{\delta}}^{(k+1)\widetilde{\delta}}m_{\delta}(u)\int_{k\widetilde{\delta}-\delta}^{u}D(D\varphi\sigma)_{i}b(X^{\delta}_{s})dsdB^{\delta,i}_{u}\\
&-\frac{2}{\gamma}\sum_{i=1}^{r}\int_{k\widetilde{\delta}}^{(k+1)\widetilde{\delta}}m_{\delta}(u)\int_{k\widetilde{\delta}-\delta}^{u}D(D\varphi\sigma)_{i}dK^{\delta}_{s}dB^{\delta,i}_{u}\\
:=&~\sum_{i=1}^{5}I_{221i}^{k}.
\end{aligned}
\end{equation}
We are going to bound each of them. By Remark $\ref{r5}$,
\begin{equation}
\begin{aligned}
E\left[\sup_{s\leq t}|\sum_{k=1}^{S_{\delta}(s)-1}I_{2215}^{k}|^{2}\right]
\lesssim~&E\left[\left(\int_{\widetilde{\delta}}^{[t](\widetilde{\delta})}|\dot{B}^{\delta}_{u}||K^{\delta}|_{[u](\widetilde{\delta})-\delta}^{[u]^{+}(\widetilde{\delta})}du\right)^{2}\right]
\lesssim E\left[\left(|K^{\delta}|_{0}^{T}\sup_{k}\int_{k\widetilde{\delta}}^{(k+1)\widetilde{\delta}}|\dot{B}^{\delta}_{u}|du\right)^{2}\right]\\
\leq~&\left(E\left[\left(|K^{\delta}|_{0}^{T}\right)^{4}\right]\right)^{\frac{1}{2}}\left(\sum_{k}E\left[\left(\int_{k\widetilde{\delta}}^{(k+1)\widetilde{\delta}}|\dot{B}^{\delta}_{u}|du\right)^{4}\right]\right)^{\frac{1}{2}}\\
\lesssim~&\left(\frac{1}{n(\delta)\delta}n(\delta)^{4}\delta^{2}\right)^{\frac{1}{2}}
\leq~n(\delta)^{\frac{3}{2}}\delta^{\frac{1}{2}}.
\end{aligned}
\end{equation}
Also
\begin{equation}
\begin{aligned}
E\left[\sup_{s\leq t}|\sum_{k=1}^{S_{\delta}(s)-1}I_{2214}^{k}|^{2}\right]
\lesssim~&E\left[\left(\int_{\widetilde{\delta}}^{[t](\widetilde{\delta})}\widetilde{\delta}|\dot{B}^{\delta}_{u}|du\right)^{2}\right]
\leq E\left[\left(\sup_{k}\int_{k\widetilde{\delta}}^{(k+1)\widetilde{\delta}}|\dot{B}^{\delta}_{u}|du\right)^{2}\right]\\
\leq~&\left(\sum_{k}E\left[\left(\int_{k\widetilde{\delta}}^{(k+1)\widetilde{\delta}}|\dot{B}^{\delta}_{u}|du\right)^{4}\right]\right)^{\frac{1}{2}}\lesssim~\left(\frac{1}{n(\delta)\delta}n(\delta)^{4}\delta^{2}\right)^{\frac{1}{2}}\\
\leq~&n(\delta)^{\frac{3}{2}}\delta^{\frac{1}{2}}.
\end{aligned}
\end{equation}
As for $I_{2212}$, we have
\begin{equation}
\begin{aligned}
E\left[\sup_{s\leq t}|\sum_{k=1}^{S_{\delta}(s)-1}I_{2212}^{k}|^{2}\right]
\lesssim&~E\left[\left(\sum_{k=1}^{S_{\delta}(t)-1}||X^{\delta}||_{[k\widetilde{\delta}-\delta,(k+1)\widetilde{\delta}]}\left(\int_{k\widetilde{\delta}-\delta}^{(k+1)\widetilde{\delta}}|\dot{B}^{\delta}_{u}|du\right)^{2}\right)^{2}\right]\\
\lesssim&~S_{\delta}(T)\sum_{k=1}^{S_{\delta}(t)-1}E\left[||X^{\delta}||_{[k\widetilde{\delta}-\delta,(k+1)\widetilde{\delta}]}^{2}\left(\int_{k\widetilde{\delta}-\delta}^{(k+1)\widetilde{\delta}}|\dot{B}^{\delta}_{u}|du\right)^{4}\right]\\
\lesssim&~S_{\delta}(T)^{2}~\widetilde{\delta}^{\frac{1}{2}}~n(\delta)^{4}\delta^{2}\lesssim n(\delta)^{\frac{5}{2}}\delta^{\frac{1}{2}}.
\end{aligned}
\end{equation}
A similar argument shows that
$$
E\left[\sup_{s\leq t}|\sum_{k=1}^{S_{\delta}(s)-1}I_{2213}^{k}|^{2}\right]\lesssim n(\delta)^{\frac{5}{2}}\delta^{\frac{1}{2}}.
$$
Observe that
\begin{equation}\label{2,2,1,1}
\begin{aligned}
&~~I_{2211}^{k}\\
=&-\frac{2}{\gamma}\sum_{i,j=1}^{r}m_{\delta}(D(D\varphi\sigma)_{i}\sigma)_{j}(X^{\delta}_{k\widetilde{\delta}-\delta})\left(B_{k\widetilde{\delta}}^{\delta,j}-B_{k\widetilde{\delta}-\delta}^{\delta,j}\right)\left(B_{(k+1)\widetilde{\delta}}^{\delta,i}-B^{\delta,i}_{k\widetilde{\delta}}\right)\\
&-\frac{2}{\gamma}\sum_{i,j=1}^{r}m_{\delta}(D(D\varphi\sigma)_{i}\sigma)_{j}(X^{\delta}_{k\widetilde{\delta}-\delta})\left(B_{(k+1)\widetilde{\delta}}^{\delta,j}-B_{k\widetilde{\delta}}^{\delta,j}\right)\left(B_{(k+1)\widetilde{\delta}}^{\delta,i}-B_{k\widetilde{\delta}}^{\delta,i}\right)\\
&+\frac{2}{\gamma}\sum_{i,j=1}^{r}\int_{k\widetilde{\delta}}^{(k+1)\widetilde{\delta}}m_{\delta}(D(D\varphi\sigma)_{i}\sigma)_{j}(X^{\delta}_{k\widetilde{\delta}-\delta})\left(B_{(k+1)\widetilde{\delta}}^{\delta,j}-B_{u}^{\delta,j}\right)dB_{u}^{\delta,i}\\
:=&\sum_{i=1}^{3}I_{2211i}^{k}
\end{aligned}
\end{equation}
Rewrite $I_{22111}^{k}$ as
\begin{equation}
\begin{aligned}
&~~~~I_{22111}^{k}\\
=&-\frac{2}{\gamma}\sum_{i,j=1}^{r}m_{\delta}(D(D\varphi\sigma)_{i}\sigma)_{j}(X^{\delta}_{k\widetilde{\delta}-\delta})\left(B_{k\widetilde{\delta}}^{\delta,j}-B_{k\widetilde{\delta}-\delta}^{\delta,j}\right)\left(W^{i}((k+1)\widetilde{\delta})-W^{i}(k\widetilde{\delta})\right)\\
&-\frac{2}{\gamma}\sum_{i,j=1}^{r}m_{\delta}(D(D\varphi\sigma)_{i}\sigma)_{j}(X^{\delta}_{k\widetilde{\delta}-\delta})\left(B_{k\widetilde{\delta}}^{\delta,j}-B_{k\widetilde{\delta}-\delta}^{\delta,j}\right)\left(B^{\delta,i}_{(k+1)\widetilde{\delta}}-W^{i}((k+1)\widetilde{\delta})\right)\\
&+\frac{2}{\gamma}\sum_{i,j=1}^{r}m_{\delta}(D(D\varphi\sigma)_{i}\sigma)_{j}(X^{\delta}_{k\widetilde{\delta}-\delta})\left(B_{k\widetilde{\delta}}^{\delta,j}-B_{k\widetilde{\delta}-\delta}^{\delta,j}\right)\left(B^{\delta,i}_{k\widetilde{\delta}}-W^{i}(k\widetilde{\delta})\right)\\
:=&\sum_{i=1}^{3}I_{22111i}^{k}.
\end{aligned}
\end{equation}
It's easy to see that
$$
E\left[\sup_{s\leq t}|\sum_{k=1}^{S_{\delta}(s)-1}\left(I_{221112}^{k}+I_{221113}^{k}\right)|^{2}\right]\lesssim n(\delta)^{-2}.
$$
Moreover, martingale inequality implies
\begin{equation}
\begin{aligned}
E\left[\sup_{s\leq t}|\sum_{k=1}^{S_{\delta}(s)-1}I_{221111}^{k}|^{2}\right]
\lesssim~&\sum_{i,j=1}^{r}\sum_{k=1}^{S_{\delta}(T)-1}E\left[\left(B^{\delta,j}_{k\widetilde{\delta}}-B^{\delta,j}_{k\widetilde{\delta}-\delta}\right)^{2}\left(W^{i}((k+1)\widetilde{\delta})-W^{i}(k\widetilde{\delta})\right)^{2}\right]\\
\lesssim~&\frac{1}{n(\delta)\delta}n(\delta)\delta^{2}=\delta.
\end{aligned}
\end{equation}
Thus
$$
E\left[\sup_{s\leq t}|\sum_{k=1}^{S_{\delta}(s)-1}I_{22111}^{k}|^{2}\right]\lesssim n(\delta)^{-2}+\delta.
$$
As for $I_{22112}^{k}$, we have
\begin{equation}
\begin{aligned}
&~~~~~~I_{22112}^{k}\\
=&-\frac{2}{\gamma}\sum_{i,j=1}^{r}m_{\delta}(D(D\varphi\sigma)_{i}\sigma)_{j}(X^{\delta}_{k\widetilde{\delta}-\delta})\left(W^{j}_{(k+1)\widetilde{\delta}-\delta}-W^{j}_{k\widetilde{\delta}-\delta}\right)
\left(W^{i}_{(k+1)\widetilde{\delta}-\delta}-W^{i}_{k\widetilde{\delta}-\delta}\right)\\
&-\frac{2}{\gamma}\sum_{i,j=1}^{r}m_{\delta}(D(D\varphi\sigma)_{i}\sigma)_{j}(X^{\delta}_{k\widetilde{\delta}-\delta})\left(W^{j}_{(k+1)\widetilde{\delta}-\delta}-W^{j}_{k\widetilde{\delta}-\delta}\right)
\left(G^{\delta,i}(0,\theta_{(k+1)\widetilde{\delta}-\delta}w)-G^{\delta,i}(0,\theta_{k\widetilde{\delta}-\delta}w)\right)\\
&-\frac{2}{\gamma}\sum_{i,j=1}^{r}m_{\delta}(D(D\varphi\sigma)_{i}\sigma)_{j}(X^{\delta}_{k\widetilde{\delta}-\delta})\left(G^{\delta,j}(0,\theta_{(k+1)\widetilde{\delta}-\delta}w)-G^{\delta,j}(0,\theta_{k\widetilde{\delta}-\delta}w)\right)
\left(W^{i}_{(k+1)\widetilde{\delta}-\delta}-W^{i}_{k\widetilde{\delta}-\delta}\right)\\
&-\frac{2}{\gamma}\sum_{i,j=1}^{r}m_{\delta}(D(D\varphi\sigma)_{i}\sigma)_{j}(X^{\delta}_{k\widetilde{\delta}-\delta})\left(G^{\delta,j}(0,\theta_{(k+1)\widetilde{\delta}-\delta}w)-G^{\delta,j}(0,\theta_{k\widetilde{\delta}-\delta}w)\right)\\
&~~~~~~~~~\times\left(G^{\delta,i}(0,\theta_{(k+1)\widetilde{\delta}-\delta}w)-G^{\delta,i}(0,\theta_{k\widetilde{\delta}-\delta}w)\right)\\
:=&~\sum_{i=1}^{4}I_{22112i}^{k}
\end{aligned}
\end{equation}
It's trivial to prove that
$$
E\left[\sup_{s\leq t}|\sum_{k=1}^{S_{\delta}(s)-1}\left(I_{221122}^{k}+I_{221123}^{k}\right)|^{2}\right]\lesssim n(\delta)^{-1}
$$
and
$$
E\left[\sup_{s\leq t}|\sum_{k=1}^{S_{\delta}(s)-1}I_{221124}^{k}|^{2}\right]\lesssim n(\delta)^{-2}.
$$
Also, by martingale inequality, we have
\begin{equation}
\begin{aligned}
E\left[\sup_{s\leq t}|\sum_{k=1}^{S_{\delta}(s)-1}I_{221121}^{k}|^{2}\right]&\lesssim n(\delta)\delta+\int_{\widetilde{\delta}}^{[t](\widetilde{\delta})}E[m_{\delta}^{2}([u](\widetilde{\delta})-\delta)]du\\
&\lesssim n(\delta)^{\frac{1}{2}}\delta^{\frac{1}{2}}+\int_{0}^{t}E[m_{\delta}^{2}(u)]du,
\end{aligned}
\end{equation}
the last $"\lesssim"$ is a result from Lemma $\ref{l8}$. Thus we say
$$
E\left[\sup_{s\leq t}|\sum_{k=1}^{S_{\delta}(s)-1}I_{22112}^{k}|^{2}\right]\lesssim n(\delta)^{\frac{1}{2}}\delta^{\frac{1}{2}}+n(\delta)^{-1}+\int_{0}^{t}E[m_{\delta}^{2}(u)]du.
$$

Now we turn to $I_{22113}^{k}$.
By Remark $\ref{r3}$, it's easy to check that
\begin{equation}
\begin{aligned}
&I_{22113}^{k}\\
=\frac{2}{\gamma}&\sum_{i,j=1}^{r}\int_{k\widetilde{\delta}-\delta}^{(k+1)\widetilde{\delta}-\delta}m_{\delta}\left(D(D\varphi\sigma)_{i}\sigma\right)_{j}(X^{\delta}_{k\widetilde{\delta}-\delta})\dot{G}_{u}^{\delta,i}\left(G_{(k+1)\widetilde{\delta}-\delta}^{\delta,j}-G_{u}^{\delta,j}\right)du\\
=\frac{2}{\gamma}&\sum_{i,j=1}^{r}\int_{k\widetilde{\delta}}^{(k+1)\widetilde{\delta}-\delta}m_{\delta}\left(D(D\varphi\sigma)_{i}\sigma\right)_{j}(X^{\delta}_{k\widetilde{\delta}-\delta})
\left[\dot{G}_{u}^{\delta,i}\left(G_{(k+1)\widetilde{\delta}-\delta}^{\delta,j}-G_{u}^{\delta,j}\right)-c_{ij}(\widetilde{\delta}-\delta,\delta)\right]du\\
+\frac{2}{\gamma}&\sum_{i,j=1}^{r}\int_{k\widetilde{\delta}-\delta}^{k\widetilde{\delta}}m_{\delta}\left(D(D\varphi\sigma)_{i}\sigma\right)_{j}(X^{\delta}_{k\widetilde{\delta}-\delta})
\left[\dot{G}_{u}^{\delta,i}\left(G_{(k+1)\widetilde{\delta}-\delta}^{\delta,j}-G_{u}^{\delta,j}\right)-c_{ij}^{*}(\widetilde{\delta},\delta)\right]du\\
+\frac{2}{\gamma}&\widetilde{\delta}\sum_{i,j=1}^{r}m_{\delta}\left(D(D\varphi\sigma)_{i}\sigma\right)_{j}(X^{\delta}_{k\widetilde{\delta}-\delta})\left(c_{ij}(\widetilde{\delta},\delta)-c_{ij}\right)\\
+\frac{2}{\gamma}&\widetilde{\delta}\sum_{i,j=1}^{r}m_{\delta}\left(D(D\varphi\sigma)_{i}\sigma\right)_{j}(X^{\delta}_{k\widetilde{\delta}-\delta})c_{ij}\\
:=\sum_{i=1}^{4}&I_{22113i}^{k}.
\end{aligned}
\end{equation}
Immediately
$$
E\left[\sup_{s\leq t}|\sum_{k=1}^{S_{\delta}(s)-1}I_{221134}^{k}|^{2}\right]
\lesssim~\int_{\widetilde{\delta}}^{[t](\widetilde{\delta})}E\left[m_{\delta}^{2}([u](\widetilde{\delta})-\delta)\right]du
\lesssim~\int_{0}^{t}E\left[m_{\delta}^{2}(u)\right]du+n(\delta)^{\frac{1}{2}}\delta^{\frac{1}{2}},
$$
the last $"\lesssim"$ is a result from Lemma $\ref{l8}$. Also we have
$$
E\left[\sup_{s\leq t}|\sum_{k=1}^{S_{\delta}(s)-1}I_{221133}^{k}|^{2}\right]
\lesssim~S_{\delta}(T)\sum_{k=1}^{S_{\delta}(t)-1}\widetilde{\delta}^{2}\sum_{i,j=1}^{r}\left(c_{ij}(\widetilde{\delta},\delta)-c_{ij}\right)^{2}
\lesssim~|(c_{ij})(\widetilde{\delta},\delta)-(c_{ij})|^{2}.
$$
Moreover
\begin{equation}
\begin{aligned}
&E\left[\sup_{s\leq t}|\sum_{k=1}^{S_{\delta}(s)-1}I_{221132}^{k}|^{2}\right]\\
\lesssim&S_{\delta}(T)\sum_{k=1}^{S_{\delta}(t)-1}\sum_{i,j=1}^{r}E\left[\left(\int_{k\widetilde{\delta}-\delta}^{k\widetilde{\delta}}\dot{G}^{\delta,i}(u,w)\left[G^{\delta,j}((k+1)\widetilde{\delta}-\delta,w)-G^{\delta,j}(u,w)\right]du-\delta c_{ij}^{*}(\widetilde{\delta},\delta)\right)^{2}\right]\\
=&S_{\delta}(T)\sum_{k=1}^{S_{\delta}(t)-1}\sum_{i,j=1}^{r}E\left[\left(\int_{0}^{\delta}\dot{G}^{\delta,i}(u,\theta_{k\widetilde{\delta}-\delta}w)\left[G^{\delta,j}(\widetilde{\delta},\theta_{k\widetilde{\delta}-\delta}w)-G^{\delta,j}(u,\theta_{k\widetilde{\delta}-\delta}w)\right]du-\delta c_{ij}^{*}(\widetilde{\delta},\delta)\right)^{2}\right]\\
\lesssim&S_{\delta}(T)^{2}\sum_{i,j=1}^{r}E\left[\left(\int_{0}^{\delta}\dot{G}^{\delta,i}(u,w)\left[G^{\delta,j}(\widetilde{\delta},w)-G^{\delta,j}(u,w)\right]du\right)^{2}\right]\\
\lesssim&S_{\delta}(T)^{2}\sum_{i,j=1}^{r}\left\{E\left[\left(\int_{0}^{\delta}\dot{G}^{\delta,i}_{u}\left[G^{\delta,j}_{\delta}-G^{\delta,j}_{u}\right]du\right)^{2}+(G^{\delta,i}_{\delta}-G^{\delta,i}_{0})^{2}(G^{\delta,j}_{\widetilde{\delta}}-G^{\delta,j}_{\delta})^{2}\right]\right\}
\\
\lesssim&S_{\delta}(T)^{2}\sum_{i,j=1}^{r}\left\{E\left[\left(\int_{0}^{\delta}~|\dot{G}^{\delta}_{u}|du\right)^{4}+\left(\int_{0}^{\delta}|\dot{G}^{\delta}_{u}|du\right)^{2}\left(G^{\delta,j}(0, \theta_{\widetilde{\delta}} w)-G^{\delta,j}(0, \theta_{\delta} w)+W(\widetilde{\delta})-W(\delta)\right)^{2}\right]\right\}\\
\lesssim&\frac{1}{n(\delta)^{2}\delta^{2}}\left(\delta^{2}+\delta^{2}n(\delta)\right)\leq~n(\delta)^{-1}.
\end{aligned}
\end{equation}
Since
$$
H_{n}=\frac{2}{\gamma}\sum_{k=1}^{n}\sum_{i,j=1}^{r}\int_{k\widetilde{\delta}}^{(k+1)\widetilde{\delta}-\delta}m_{\delta}\left(D(D\varphi\sigma)_{i}\sigma\right)_{j}(X^{\delta}_{k\widetilde{\delta}-\delta})
\left[\dot{G}_{u}^{\delta,i}\left(G_{(k+1)\widetilde{\delta}-\delta}^{\delta,j}-G_{u}^{\delta,j}\right)-c_{ij}(\widetilde{\delta}-\delta,\delta)\right]du
$$
is a $\mathcal{H}_{n}-$martingale, where $\mathcal{H}_{n}=\mathcal{B}_{(n+1)\widetilde{\delta}}$, by martingale inequality, we have
\begin{equation}
\begin{aligned}
&E\left[\sup_{s\leq t}|\sum_{k=1}^{S_{\delta}(s)-1}I_{221131}^{k}|^{2}\right]\\
\lesssim~&\sum_{i,j=1}^{r}\sum_{k=1}^{S_{\delta}(T)-1}E\left[\left(\int_{k\widetilde{\delta}}^{(k+1)\widetilde{\delta}-\delta}m_{\delta}\left(D(D\varphi\sigma)_{i}\sigma\right)_{j}(X^{\delta}_{k\widetilde{\delta}-\delta})\left[\dot{G}_{u}^{\delta,i}\left(G_{(k+1)\widetilde{\delta}-\delta}^{\delta,j}-G_{u}^{\delta,j}\right)-c_{ij}(\widetilde{\delta}-\delta,\delta)\right]du\right)^{2}\right]\\
\lesssim~&S_{\delta}(T)\sum_{i,j=1}^{r}E\left[\left(\int_{0}^{\widetilde{\delta}-\delta}\left[\dot{G}_{u}^{\delta,i}\left(G_{\widetilde{\delta}-\delta}^{\delta,j}-G_{u}^{\delta,j}\right)-c_{ij}(\widetilde{\delta}-\delta,\delta)\right]du\right)^{2}\right]\\
\lesssim~&S_{\delta}(T)\sum_{i,j=1}^{r}E\left[\left(\int_{0}^{\widetilde{\delta}-\delta}\dot{G}_{u}^{\delta,i}\left[G_{\widetilde{\delta}-\delta}^{\delta,j}-G_{u}^{\delta,j}\right]du\right)^{2}\right]\\
\lesssim~&S_{\delta}(T)E\left[\left(\int_{0}^{\widetilde{\delta}}|\dot{G}^{\delta}_{u}|du\right)^{4}\right]\lesssim n(\delta)^{3}\delta.
\end{aligned}
\end{equation}
That is
$$
E\left[\sup_{s\leq t}|\sum_{k=1}^{S_{\delta}(s)-1}I_{22113}^{k}|^{2}\right]
\lesssim\int_{0}^{t}E\left[m_{\delta}^{2}(u)\right]du+o(1).
$$
Putting our estimates together, we have
$$
E\left[\sup_{s\leq t}|\sum_{k=1}^{S_{\delta}(s)-1}I_{221}^{k}|^{2}\right]
\lesssim\int_{0}^{t}E\left[m_{\delta}^{2}(u)\right]du+o(1).
$$

Now we consider $I_{222}^{k}$.
\begin{equation}
\begin{aligned}
&~~I_{222}^{k}\\
=&-\frac{2}{\gamma}\sum_{i=1}^{r}\int_{k\widetilde{\delta}}^{(k+1)\widetilde{\delta}}(D\varphi\sigma)_{i}(X^{\delta}_{k\widetilde{\delta}-\delta})\left(m_{\delta}(u)-m_{\delta}(k\widetilde{\delta}-\delta)\right)dB_{u}^{\delta,i}\\
=&\frac{4}{\gamma^{2}}\sum_{i=1}^{r}\int_{k\widetilde{\delta}}^{(k+1)\widetilde{\delta}}(D\varphi\sigma)_{i}(X^{\delta}_{k\widetilde{\delta}-\delta})\int_{k\widetilde{\delta}-\delta}^{u}m_{\delta}(s)I_{\delta}(s)dsdB_{u}^{\delta,i}\\
+&\frac{4}{\gamma^{2}}\sum_{i=1}^{r}\int_{k\widetilde{\delta}}^{(k+1)\widetilde{\delta}}(D\varphi\sigma)_{i}(X^{\delta}_{k\widetilde{\delta}-\delta})\int_{k\widetilde{\delta}-\delta}^{u}m_{\delta}(s)D\varphi\sigma(X^{\delta}_{s})dB^{\delta}_{s}dB_{u}^{\delta,i}\\
+&\frac{4}{\gamma^{2}}\sum_{i=1}^{r}\int_{k\widetilde{\delta}}^{(k+1)\widetilde{\delta}}(D\varphi\sigma)_{i}(X^{\delta}_{k\widetilde{\delta}-\delta})\int_{k\widetilde{\delta}-\delta}^{u}m_{\delta}(s)D\varphi\sigma(X_{s})dW_{s}dB_{u}^{\delta,i}
\\
-&\frac{2}{\gamma}\sum_{i=1}^{r}\int_{k\widetilde{\delta}}^{(k+1)\widetilde{\delta}}(D\varphi\sigma)_{i}(X^{\delta}_{k\widetilde{\delta}-\delta})\int_{k\widetilde{\delta}-\delta}^{u}\mu_{\delta}(s)dA^{\delta}_{s}dB_{u}^{\delta,i}\\
-&\frac{4}{\gamma}\sum_{i=1}^{r}\int_{k\widetilde{\delta}}^{(k+1)\widetilde{\delta}}(D\varphi\sigma)_{i}(X^{\delta}_{k\widetilde{\delta}-\delta})\int_{k\widetilde{\delta}-\delta}^{u}\mu_{\delta}(s)(X_{s}-X^{\delta}_{s})^{*}(b(X_{s})-b(X^{\delta}_{s}))dsdB_{u}^{\delta,i}\\
+&\frac{8}{\gamma^{2}}\sum_{i=1}^{r}\int_{k\widetilde{\delta}}^{(k+1)\widetilde{\delta}}(D\varphi\sigma)_{i}(X^{\delta}_{k\widetilde{\delta}-\delta})\int_{k\widetilde{\delta}-\delta}^{u}\mu_{\delta}(s)D\varphi\sigma\sigma^{*}(X_{s})(X_{s}-X^{\delta}_{s})dsdB_{u}^{\delta,i}\\
+&\frac{4}{\gamma}\sum_{i=1}^{r}\int_{k\widetilde{\delta}}^{(k+1)\widetilde{\delta}}(D\varphi\sigma)_{i}(X^{\delta}_{k\widetilde{\delta}-\delta})\int_{k\widetilde{\delta}-\delta}^{u}\mu_{\delta}(s)(X_{s}-X^{\delta}_{s})^{*}\sigma(X^{\delta}_{s})dB^{\delta}_{s}dB_{u}^{\delta,i}
\end{aligned}
\end{equation}	
\begin{equation}
\begin{aligned}
-&\frac{4}{\gamma}\sum_{i=1}^{r}\int_{k\widetilde{\delta}}^{(k+1)\widetilde{\delta}}(D\varphi\sigma)_{i}(X^{\delta}_{k\widetilde{\delta}-\delta})\int_{k\widetilde{\delta}-\delta}^{u}\mu_{\delta}(s)(X_{s}-X^{\delta}_{s})^{*}\sigma(X_{s})dW_sdB_{u}^{\delta,i}\\
-&\frac{4}{\gamma}\sum_{i=1}^{r}\int_{k\widetilde{\delta}}^{(k+1)\widetilde{\delta}}(D\varphi\sigma)_{i}(X^{\delta}_{k\widetilde{\delta}-\delta})\int_{k\widetilde{\delta}-\delta}^{u}\mu_{\delta}(s)(X_{s}-X^{\delta}_{s})^{*}\left(\bar{b}-b\right)(X_{s})dsdB_{u}^{\delta,i}\\
-&\frac{2}{\gamma}\sum_{i=1}^{r}\int_{k\widetilde{\delta}}^{(k+1)\widetilde{\delta}}(D\varphi\sigma)_{i}(X^{\delta}_{k\widetilde{\delta}-\delta})\int_{k\widetilde{\delta}-\delta}^{u}\mu_{\delta}(s)tr(\sigma^{*}\sigma(X_{s}))dsdB_{u}^{\delta,i}\\
:=&\sum_{i=1}^{10}I_{222i}^{k},
\end{aligned}
\end{equation}
where
$$
I_{\delta}(s)=D\varphi(X^{\delta}_{s})b(X^{\delta}_{s})+D\varphi(X_{s})\bar{b}(X_{s})
+\frac{1}{2}tr(\sigma^{*}D^{2}\varphi\sigma(X_{s}))-\frac{1}{\gamma}|D\varphi(X_{s})b(X_{s})|^{2},
$$
and
$$
dA^{\delta}_{s}=-\frac{2}{\gamma}|X_{s}-X^{\delta}_{s}|^{2}(D\varphi(X^{\delta}_{s})dK^{\delta}_{s}+ D\varphi(X_{s})dK_{s})
+2( X_{s}-X^{\delta}_{s})^{*}(dK_{s}-dK^{\delta}_{s}).
$$
By H\"older's inequality, we have
$$
E\left[\sup_{s\leq t}|\sum_{k=1}^{S_{\delta}(s)-1}I_{222i}^{k}|^{2}\right]\lesssim E\left[\sup_{k}\left(\int_{k\widetilde{\delta}-\delta}^{(k+1)\widetilde{\delta}}|\dot{B}^{\delta}_{u}|du\right)^{2}\right]
\lesssim\left(\sum_{k}E\left[\left(\int_{k\widetilde{\delta}-\delta}^{(k+1)\widetilde{\delta}}|\dot{B}^{\delta}_{u}|du\right)^{4}\right]\right)^{\frac{1}{2}}\leq n(\delta)^{\frac{3}{2}}\delta^{\frac{1}{2}},
$$
for $i=1,5,6,9,10$. Also it's easy to see that
$$
E\left[\sup_{s\leq t}|\sum_{k=1}^{S_{\delta}(s)-1}I_{2224}^{k}|^{2}\right]\lesssim E\left[\left(\sup_{k}\int_{k\widetilde{\delta}}^{(k+1)\widetilde{\delta}}|_{\delta}\dot{B}_{u}|du|_{\delta}A|_{0}^{T}\right)^{2}\right]
\lesssim\left(\sum_{k}E\left[\left(\int_{k\widetilde{\delta}}^{(k+1)\widetilde{\delta}}|_{\delta}\dot{B}_{u}|du\right)^{4}\right]\right)^{\frac{1}{2}}\leq n(\delta)^{\frac{3}{2}}\delta^{\frac{1}{2}}.
$$

Note that $I_{2223}^{k}$ equals to
\begin{equation}\label{2,2,2,3}
\begin{aligned}
&\frac{4}{\gamma^{2}}\sum_{i=1}^{r}\int_{k\widetilde{\delta}}^{(k+1)\widetilde{\delta}}(D\varphi\sigma)_{i}(X^{\delta}_{k\widetilde{\delta}-\delta})\int_{k\widetilde{\delta}-\delta}^{u}m_{\delta}D\varphi\sigma(X_{s})dW_{s}dB_{u}^{\delta,i}\\
=&\frac{4}{\gamma^{2}}\sum_{i=1}^{r}\int_{k\widetilde{\delta}}^{(k+1)\widetilde{\delta}}(D\varphi\sigma)_{i}(X^{\delta}_{k\widetilde{\delta}-\delta})\int_{k\widetilde{\delta}-\delta}^{u}m_{\delta}D\varphi\sigma(X_{s})-m_{\delta}D\varphi\sigma(X_{k\widetilde{\delta}-\delta})dW_{s}dB_{u}^{\delta,i}
\\
+&\frac{4}{\gamma^{2}}\sum_{i,j=1}^{r}m_{\delta}(D\varphi\sigma)_{i}(X^{\delta}_{k\widetilde{\delta}-\delta})(D\varphi\sigma)_{j}(X_{k\widetilde{\delta}-\delta})\left(W_{(k+1)\widetilde{\delta}}^{j}-W_{k\widetilde{\delta}-\delta}^{j}\right)\left(B_{(k+1)\widetilde{\delta}}^{\delta,i}-B_{k\widetilde{\delta}}^{\delta,i}\right)\\
-&\frac{4}{\gamma^{2}}\sum_{i,j=1}^{r}m_{\delta}(D\varphi\sigma)_{i}(X^{\delta}_{k\widetilde{\delta}-\delta})(D\varphi\sigma)_{j}(X_{k\widetilde{\delta}-\delta})\int_{k\widetilde{\delta}}^{(k+1)\widetilde{\delta}}\left(B^{\delta,i}_{u}-B^{\delta,i}_{k\widetilde{\delta}}\right)dW_{u}^{j}\\
:=&\sum_{i=1}^{3}I_{2223i}^{k}.
\end{aligned}
\end{equation}
H\"older's inequality and BDG's inequality yield
\begin{equation}
\begin{aligned}
&~~E\left[\sup_{s\leq t}|\sum_{k=1}^{S_{\delta}(s)-1}I_{22231}^{k}|^{2}\right]
\\
\lesssim~&S_{\delta}(T)\sum_{k=1}^{S_{\delta}(T)-1}E\left[\left(\int_{k\widetilde{\delta}}^{(k+1)\widetilde{\delta}}|\dot{B}_{u}^{\delta,i}|du\right)^{2}
\sup_{k\widetilde{\delta}\leq
u\leq(k+1)\widetilde{\delta}}|\int_{k\widetilde{\delta}-\delta}^{u}m_{\delta}D\varphi\sigma(X_{s})-m_{\delta}D\varphi\sigma(X_{k\widetilde{\delta}-\delta})dW_{s}|^{2}\right]\\
\lesssim~&S_{\delta}(T)^{2}n(\delta)^{2}\delta\widetilde{\delta}^{\frac{3}{2}}\lesssim~ n(\delta)^{\frac{3}{2}}\delta^{\frac{1}{2}}
\end{aligned}
\end{equation}
and
\begin{equation}
\begin{aligned}
&E\left[\sup_{s\leq t}|\sum_{k=1}^{S_{\delta}(s)-1}I_{22233}^{k}|^{2}\right]\\
\lesssim~&\sum_{i,j=1}^{r}E\left[\sup_{\widetilde{\delta}\leq s\leq t}|\int_{\widetilde{\delta}}^{s}m_{\delta}(D\varphi\sigma)_{i}(X^{\delta}_{[u](\widetilde{\delta})-\delta})(D\varphi\sigma)_{j}(X_{[u](\widetilde{\delta})-\delta})\left(B^{\delta,i}_{u}-B^{\delta,i}_{[u](\widetilde{\delta})}\right)dW_{u}^{j}|^{2}\right]\\
\lesssim~&\int_{0}^{T}E\left[\left(\int_{[u](\widetilde{\delta})}^{[u]^{+}(\widetilde{\delta})}|\dot{B}^{\delta}_{s}|ds\right)^{2}\right]du\lesssim n(\delta)^{2}\delta.
\end{aligned}
\end{equation}
We split $I_{22232}^{k}$ as
\begin{equation}
\begin{aligned}
&I_{22232}^{k}\\
=~&\frac{4}{\gamma^{2}}\sum_{i,j=1}^{r}m_{\delta}(D\varphi\sigma)_{i}(X^{\delta}_{k\widetilde{\delta}-\delta})(D\varphi\sigma)_{j}(X_{k\widetilde{\delta}-\delta})\left(W_{(k+1)\widetilde{\delta}-\delta}^{j}-W_{k\widetilde{\delta}-\delta}^{j}\right)\left(W_{(k+1)\widetilde{\delta}-\delta}^{i}-W_{k\widetilde{\delta}-\delta}^{i}\right)
\\
+&\frac{4}{\gamma^{2}}\sum_{i,j=1}^{r}m_{\delta}(D\varphi\sigma)_{i}(X^{\delta}_{k\widetilde{\delta}-\delta})(D\varphi\sigma)_{j}(X_{k\widetilde{\delta}-\delta})\left(W_{(k+1)\widetilde{\delta}-\delta}^{j}-W_{k\widetilde{\delta}-\delta}^{j}\right)\left(G^{\delta,i}(0,\theta_{(k+1)\widetilde{\delta}-\delta}w)-G^{\delta,i}(0,\theta_{k\widetilde{\delta}-\delta}w)\right)
\\
+&\frac{4}{\gamma^{2}}\sum_{i,j=1}^{r}m_{\delta}(D\varphi\sigma)_{i}(X^{\delta}_{k\widetilde{\delta}-\delta})(D\varphi\sigma)_{j}(X_{k\widetilde{\delta}-\delta})\left(W_{(k+1)\widetilde{\delta}}^{j}-W_{(k+1)\widetilde{\delta}-\delta}^{j}\right)\left(W_{(k+1)\widetilde{\delta}-\delta}^{i}-W_{k\widetilde{\delta}-\delta}^{i}\right)
\\
+&\frac{4}{\gamma^{2}}\sum_{i,j=1}^{r}m_{\delta}(D\varphi\sigma)_{i}(X^{\delta}_{k\widetilde{\delta}-\delta})(D\varphi\sigma)_{j}(X_{k\widetilde{\delta}-\delta})\left(W_{(k+1)\widetilde{\delta}}^{j}-W_{(k+1)\widetilde{\delta}-\delta}^{j}\right)\left(G^{\delta,i}(0,\theta_{(k+1)\widetilde{\delta}-\delta}w)-G^{\delta,i}(0,\theta_{k\widetilde{\delta}-\delta}w)\right)
\\
:=~&\sum_{i=1}^{4}I_{22232i}^{k},
\end{aligned}
\end{equation}
reasoning in the same way as we handle $I_{22112}^{k}$, it's trivial to prove that
$$
E\left[\sup_{s\leq t}|\sum_{k=1}^{S_{\delta}(s)-1}I_{22232}^{k}|^{2}\right]\lesssim n(\delta)^{\frac{1}{2}}\delta^{\frac{1}{2}}+n(\delta)^{-1}+\int_{0}^{t}E[m_{\delta}^{2}(u)]du.
$$
Thus we claim
$$
E\left[\sup_{s\leq t}|\sum_{k=1}^{S_{\delta}(s)-1}I_{2223}^{k}|^{2}\right]
\lesssim\int_{0}^{t}E\left[m_{\delta}^{2}(u)\right]du+o(1).
$$
A similar argument can change $I_{2228}^{k}$ into
\begin{multline}
I_{2228}^{k,1}:=-\frac{4}{\gamma}\sum_{i,j=1}^{r}\mu_{\delta}(D\varphi\sigma)_{i}(X^{\delta}_{k\widetilde{\delta}-\delta})((X_{k\widetilde{\delta}-\delta}-X^{\delta}_{k\widetilde{\delta}-\delta})^{*}\sigma)_{j}(X_{k\widetilde{\delta}-\delta})\\
\times(W^{i}_{(k+1)\widetilde{\delta}-\delta}-W^{i}_{k\widetilde{\delta}-\delta})(W^{j}_{(k+1)\widetilde{\delta}-\delta}-W^{j}_{k\widetilde{\delta}-\delta}),
\end{multline}
if we neglect some trivial terms.

Applying the same way that we use to deal with $I_{2211}^{k}$, we can prove
$$
E\left[\sup_{s\leq t}|\sum_{k=1}^{S_{\delta}(s)-1}I_{2222}^{k}|^{2}\right]
\lesssim\int_{0}^{t}E\left[m_{\delta}^{2}(u)\right]du+o(1),
$$
and turn $I_{2227}^{k}$ into
\begin{equation}
\begin{aligned}
&I_{2227}^{k,2}+I_{2227}^{k,3}\\
=-&\frac{4}{\gamma}\sum_{i,j=1}^{r}\mu_{\delta}(D\varphi\sigma)_{i}((X_{k\widetilde{\delta}-\delta}-X^{\delta}_{k\widetilde{\delta}-\delta})^{*}\sigma)_{j}(X^{\delta}_{k\widetilde{\delta}-\delta})c_{ij}\widetilde{\delta}\\
+&\frac{4}{\gamma}\sum_{i,j=1}^{r}\mu_{\delta}(D\varphi\sigma)_{i}((X_{k\widetilde{\delta}-\delta}-X^{\delta}_{k\widetilde{\delta}-\delta})^{*}\sigma)_{j}(X^{\delta}_{k\widetilde{\delta}-\delta})(W^{i}_{(k+1)\widetilde{\delta}-\delta}-W^{i}_{k\widetilde{\delta}-\delta})(W^{j}_{(k+1)\widetilde{\delta}-\delta}-W^{j}_{k\widetilde{\delta}-\delta}).
\end{aligned}
\end{equation}

Now we consider $I_{223}^{k}$, which is the last term in $I_{22}^{k}$.
\begin{equation}
\begin{aligned}
I_{223}^{k}=&-\frac{2}{\gamma}\sum_{i=1}^{r}m_{\delta}(D\varphi\sigma)_{i}(X^{\delta}_{k\widetilde{\delta}-\delta})(W_{(k+1)\widetilde{\delta}-\delta}^{i}-W_{k\widetilde{\delta}-\delta}^{i})\\
&+\frac{2}{\gamma}\sum_{i=1}^{r}m_{\delta}(D\varphi\sigma)_{i}(X^{\delta}_{k\widetilde{\delta}-\delta})(G^{\delta,i}(0,\theta_{(k+1)\widetilde{\delta}-\delta}w)-G^{\delta,i}(0,\theta_{k\widetilde{\delta}-\delta}w))\\
:=&\sum_{i=1}^{2}I_{223i}^{k}.
\end{aligned}
\end{equation}
BDG's inequality shows
$$
E\left[\sup_{s\leq t}|\sum_{k=1}^{S_{\delta}(s)-1}I_{2231}^{k}|^{2}\right]
\lesssim\int_{0}^{t}E\left[m_{\delta}^{2}(u)\right]du+o(1),
$$
and martingale inequality gives
$$
E\left[\sup_{s\leq t}|\sum_{k=1}^{S_{\delta}(s)-1}I_{2232}^{k}|^{2}\right]\lesssim n(\delta)^{-1}.
$$

Write $I_{7}(t)$ as
\begin{equation}
\begin{aligned}
&~I_{7}(t)\\
=&-2\int_{[t](\widetilde{\delta})}^{t}\mu_{\delta}(u)(X_{u}-X^{\delta}_{u})^{*}\sigma(X^{\delta}_{u})dB^{\delta}_{u}\\
&-2\int_{\widetilde{\delta}}^{[t](\widetilde{\delta})}\mu_{\delta}(u)(X_{u}-X^{\delta}_{u})^{*}\sigma(X^{\delta}_{u})dB^{\delta}_{u}\\
&-2\int_{0}^{\widetilde{\delta}}\mu_{\delta}(u)(X_{u}-X^{\delta}_{u})^{*}\sigma(X^{\delta}_{u})dB^{\delta}_{u}\\
:=&~\sum_{i=1}^{3}I_{7i}(t).
\end{aligned}
\end{equation}
Obviously
$$
E\left[\sup_{s\leq t}|I_{71}(s)+I_{73}(s)|^{2}\right]\leq\left(\sum_{k}E\left[\left(\int_{k\widetilde{\delta}}^{(k+1)\widetilde{\delta}}|\dot{B}^{\delta}_{u}|du\right)^{4}\right]\right)^{\frac{1}{2}}\lesssim n(\delta)^{\frac{3}{2}}\delta^{\frac{1}{2}}
$$
and
\begin{equation}
\begin{aligned}
&~I_{72}(t)\\
=&-2\sum_{k=1}^{S_{\delta}(t)-1}\int_{k\widetilde{\delta}}^{(k+1)\widetilde{\delta}}(\mu_{\delta}(u)-\mu_{\delta}(k\widetilde{\delta}-\delta))(X_{u}-X^{\delta}_{u})^{*}\sigma(X^{\delta}_{u})dB^{\delta}_{u}\\
&-2\sum_{k=1}^{S_{\delta}(t)-1}\int_{k\widetilde{\delta}}^{(k+1)\widetilde{\delta}}\mu_{\delta}(k\widetilde{\delta}-\delta)\left[(X_{u}-X^{\delta}_{u})^{*}-(X_{k\widetilde{\delta}-\delta}-X^{\delta}_{k\widetilde{\delta}-\delta})^{*}\right]\sigma(X^{\delta}_{u})dB^{\delta}_{u}\\
&-2\sum_{k=1}^{S_{\delta}(t)-1}\int_{k\widetilde{\delta}}^{(k+1)\widetilde{\delta}}\mu_{\delta}(k\widetilde{\delta}-\delta)(X_{k\widetilde{\delta}-\delta}-X^{\delta}_{k\widetilde{\delta}-\delta})^{*}\left(\sigma(X^{\delta}_{u})-\sigma(X^{\delta}_{k\widetilde{\delta}-\delta})\right)dB^{\delta}_{u}\\
&-2\sum_{k=1}^{S_{\delta}(t)-1}\int_{k\widetilde{\delta}}^{(k+1)\widetilde{\delta}}\mu_{\delta}(k\widetilde{\delta}-\delta)(X_{k\widetilde{\delta}-\delta}-X^{\delta}_{k\widetilde{\delta}-\delta})^{*}\sigma(X^{\delta}_{k\widetilde{\delta}-\delta})dB^{\delta}_{u}\\
:=&~\sum_{i=1}^{4}I_{72i}^{k}.
\end{aligned}
\end{equation}
We are going to bound each of them.

With respect to $I_{724}^{k}$, write
\begin{equation}
\begin{aligned}
&~I_{724}^{k}\\
=&-2\mu_{\delta}(k\widetilde{\delta}-\delta)(X_{k\widetilde{\delta}-\delta}-X^{\delta}_{k\widetilde{\delta}-\delta})^{*}\sigma(X^{\delta}_{k\widetilde{\delta}-\delta})\left(W_{(k+1)\widetilde{\delta}-\delta}-W_{k\widetilde{\delta}-\delta}\right)\\
&-2\mu_{\delta}(k\widetilde{\delta}-\delta)(X_{k\widetilde{\delta}-\delta}-X^{\delta}_{k\widetilde{\delta}-\delta})^{*}\sigma(X^{\delta}_{k\widetilde{\delta}-\delta})\left(G^{\delta}(0,\theta_{(k+1)\widetilde{\delta}-\delta}w)-G^{\delta}(0,\theta_{k\widetilde{\delta}-\delta}w)\right)\\
:=&~I_{7241}^{4,k}+I_{7242}^{k}.
\end{aligned}
\end{equation}
Martingale inequality yields
$$
E\left[\sup_{s\leq t}|\sum_{k=1}^{S_{\delta}(s)-1}I_{7242}^{k}|^{2}\right]\lesssim n(\delta)^{-1}.
$$

Rewrite $I_{722}^{k}$ as
\begin{equation}
\begin{aligned}
&~I_{722}^{k}\\
=&-2\sum_{k=1}^{S_{\delta}(t)-1}\int_{k\widetilde{\delta}}^{(k+1)\widetilde{\delta}}\mu_{\delta}(k\widetilde{\delta}-\delta)\left[\int_{k\widetilde{\delta}-\delta}^{u}\sigma(X_{s})dW_{s}\right]^{*}\sigma(X^{\delta}_{u})dB^{\delta}_{u}\\
&-2\sum_{k=1}^{S_{\delta}(t)-1}\int_{k\widetilde{\delta}}^{(k+1)\widetilde{\delta}}\mu_{\delta}(k\widetilde{\delta}-\delta)\left[\int_{k\widetilde{\delta}-\delta}^{u}\bar{b}(X_{s})ds\right]^{*}\sigma(X^{\delta}_{u})dB^{\delta}_{u}\\
&-2\sum_{k=1}^{S_{\delta}(t)-1}\int_{k\widetilde{\delta}}^{(k+1)\widetilde{\delta}}\mu_{\delta}(k\widetilde{\delta}-\delta)\left[\int_{k\widetilde{\delta}-\delta}^{u}dK_{s}\right]^{*}\sigma(X^{\delta}_{u})dB^{\delta}_{u}\\
&+2\sum_{k=1}^{S_{\delta}(t)-1}\int_{k\widetilde{\delta}}^{(k+1)\widetilde{\delta}}\mu_{\delta}(k\widetilde{\delta}-\delta)\left[\int_{k\widetilde{\delta}-\delta}^{u}\sigma(X^{\delta}_{s})dB^{\delta}_{s}\right]^{*}\sigma(X^{\delta}_{u})dB^{\delta}_{u}\\
&+2\sum_{k=1}^{S_{\delta}(t)-1}\int_{k\widetilde{\delta}}^{(k+1)\widetilde{\delta}}\mu_{\delta}(k\widetilde{\delta}-\delta)\left[\int_{k\widetilde{\delta}-\delta}^{u}b(X^{\delta}_{s})ds\right]^{*}\sigma(X^{\delta}_{u})dB^{\delta}_{u}\\
&+2\sum_{k=1}^{S_{\delta}(t)-1}\int_{k\widetilde{\delta}}^{(k+1)\widetilde{\delta}}\mu_{\delta}(k\widetilde{\delta}-\delta)\left[\int_{k\widetilde{\delta}-\delta}^{u}dK^{\delta}_{s}\right]^{*}\sigma(X^{\delta}_{u})dB^{\delta}_{u}\\
:=&~\sum_{i=1}^{6}I_{722i}^{k}.
\end{aligned}
\end{equation}
Immediately, we have
$$
E\left[\sup_{s\leq t}|\sum_{k=1}^{S_{\delta}(s)-1}(I_{7222}^{k}+I_{7225}^{k})|^{2}\right]
\lesssim E\left[\left(\int_{\widetilde{\delta}}^{[t](\widetilde{\delta})}\widetilde{\delta}|\dot{B}^{\delta}_{u}|du\right)^{2}\right]
\lesssim~
\left(\sum_{k}E\left[\left(\int_{k\widetilde{\delta}}^{(k+1)\widetilde{\delta}}|\dot{B}^{\delta}_{u}|du\right)^{4}\right]\right)^{\frac{1}{2}}\lesssim n(\delta)^{\frac{3}{2}}\delta^{\frac{1}{2}}
$$
and
\begin{equation}
\begin{aligned}
E\left[\sup_{s\leq t}|\sum_{k=1}^{S_{\delta}(s)-1}(I_{7223}^{k}+I_{7226}^{k})|^{2}\right]
\lesssim&~ E\left[\left(|V^{\delta}|_{0}^{T}\sup_{k}\int_{k\widetilde{\delta}}^{(k+1)\widetilde{\delta}}|\dot{B}^{\delta}_{u}|du\right)^{2}\right]\\
\lesssim&~
\left(\sum_{k}E\left[\left(\int_{k\widetilde{\delta}}^{(k+1)\widetilde{\delta}}|\dot{B}^{\delta}_{u}|du\right)^{4}\right]\right)^{\frac{1}{2}}\lesssim n(\delta)^{\frac{3}{2}}\delta^{\frac{1}{2}},
\end{aligned}
\end{equation}
where $V^{\delta}_{t}:=~K^{\delta}_{t}+K_{t}$.

Employing the same way that we deal with $I_{2223}^{k}$ in $(\ref{2,2,2,3})$, turn $I_{7221}^{k}$ into
$$
I_{7221}^{5,k}:=-2\mu_{\delta}(k\widetilde{\delta}-\delta)\left(W_{(k+1)\widetilde{\delta}-\delta}-W_{k\widetilde{\delta}-\delta}\right)^{*}\sigma^{*}(X_{k\widetilde{\delta}-\delta})\sigma(X^{\delta}_{k\widetilde{\delta}-\delta})\left(W_{(k+1)\widetilde{\delta}-\delta}-W_{k\widetilde{\delta}-\delta}\right).
$$
Using the way in $(\ref{2,2,1,1})$ to handle $I_{7224}^{k}$, we can change $I_{7224}^{k}$ into
\begin{equation}
\begin{aligned}
&2\mu_{\delta}(k\widetilde{\delta}-\delta)\left(W_{(k+1)\widetilde{\delta}-\delta}-W_{k\widetilde{\delta}-\delta}\right)^{*}\sigma^{*}\sigma(X^{\delta}_{k\widetilde{\delta}-\delta})\left(W_{(k+1)\widetilde{\delta}-\delta}-W_{k\widetilde{\delta}-\delta}\right)\\
-&2\mu_{\delta}(k\widetilde{\delta}-\delta)\sum_{i,j=1}^{r}\left(\sigma^{*}\sigma\right)^{i}_{j}(X^{\delta}_{k\widetilde{\delta}-\delta})c_{ij}\widetilde{\delta}\\
:=&I_{72241}^{6,k}+I_{72242}^{7,k}
\end{aligned}
\end{equation}
without some trivial terms.

Now we apply It\^o's formula to $I_{721}^{k}$, it follows that
\begin{equation}
\begin{aligned}
&-2\int_{k\widetilde{\delta}}^{(k+1)\widetilde{\delta}}(\mu_{\delta}(u)-\mu_{\delta}(k\widetilde{\delta}-\delta))(X_{u}-X^{\delta}_{u})^{*}\sigma(X^{\delta}_{u})dB^{\delta}_{u}\\
=&\frac{4}{\gamma}\int_{k\widetilde{\delta}}^{(k+1)\widetilde{\delta}}(\int_{k\widetilde{\delta}-\delta}^{u}\mu_{\delta}D\varphi b(X^{\delta}_{s})ds)(X_{u}-X^{\delta}_{u})^{*}\sigma(X^{\delta}_{u})dB^{\delta}_{u}\\
+&\frac{4}{\gamma}\int_{k\widetilde{\delta}}^{(k+1)\widetilde{\delta}}(\int_{k\widetilde{\delta}-\delta}^{u}\mu_{\delta}D\varphi\sigma(X^{\delta}_{s})dB^{\delta}_{s})(X_{u}-X^{\delta}_{u})^{*}\sigma(X^{\delta}_{u})dB^{\delta}_{u}\\
+&\frac{4}{\gamma}\int_{k\widetilde{\delta}}^{(k+1)\widetilde{\delta}}(\int_{k\widetilde{\delta}-\delta}^{u}\mu_{\delta}D\varphi(X^{\delta}_{s})dK^{\delta}_{s})(X_{u}-X^{\delta}_{u})^{*}\sigma(X^{\delta}_{u})dB^{\delta}_{u}\\
+&\frac{4}{\gamma}\int_{k\widetilde{\delta}}^{(k+1)\widetilde{\delta}}(\int_{k\widetilde{\delta}-\delta}^{u}\mu_{\delta}D\varphi \bar{b}(X_{s})ds)(X_{u}-X^{\delta}_{u})^{*}\sigma(X^{\delta}_{u})dB^{\delta}_{u}\\
+&\frac{4}{\gamma}\int_{k\widetilde{\delta}}^{(k+1)\widetilde{\delta}}(\int_{k\widetilde{\delta}-\delta}^{u}\mu_{\delta}D\varphi\sigma(X_{s})dW_{s})(X_{u}-X^{\delta}_{u})^{*}\sigma(X^{\delta}_{u})dB^{\delta}_{u}\\
+&\frac{4}{\gamma}\int_{k\widetilde{\delta}}^{(k+1)\widetilde{\delta}}(\int_{k\widetilde{\delta}-\delta}^{u}\mu_{\delta}D\varphi(X_{s})dK_{s})(X_{u}-X^{\delta}_{u})^{*}\sigma(X^{\delta}_{u})dB^{\delta}_{u}\\
+&\frac{2}{\gamma}\int_{k\widetilde{\delta}}^{(k+1)\widetilde{\delta}}(\int_{k\widetilde{\delta}-\delta}^{u}\mu_{\delta}tr(\sigma^{*}D^{2}\varphi\sigma(X_{s}))ds)(X_{u}-X^{\delta}_{u})^{*}\sigma(X^{\delta}_{u})dB^{\delta}_{u}\\
-&\frac{4}{\gamma^{2}}\int_{k\widetilde{\delta}}^{(k+1)\widetilde{\delta}}(\int_{k\widetilde{\delta}-\delta}^{u}\mu_{\delta}|D\varphi\sigma(X_{s})|^{2}ds)(X_{u}-X^{\delta}_{u})^{*}\sigma(X^{\delta}_{u})dB^{\delta}_{u}
\\
:=&~\sum_{i=1}^{8}I_{721i}^{k}.
\end{aligned}
\end{equation}
Then we obtain
$$
E\left[\sup_{s\leq t}|\sum_{k=1}^{S_{\delta}(s)-1}I_{721i}^{k}|^{2}\right]
\lesssim~
\left(\sum_{k}E\left[\left(\int_{k\widetilde{\delta}}^{(k+1)\widetilde{\delta}}|\dot{B}^{\delta}_{u}|du\right)^{4}\right]\right)^{\frac{1}{2}}\lesssim n(\delta)^{\frac{3}{2}}\delta^{\frac{1}{2}},
$$
for $i=1,3,4,6,7,8$.
As we treat $I_{2211}^{k}$ in $(\ref{2,2,1,1})$, $I_{7212}^{k}$ can be considered as
\begin{equation}
\begin{aligned}
&\frac{4}{\gamma}\sum_{i,j=1}^{r}\mu_{\delta}(k\widetilde{\delta}-\delta)(D\varphi\sigma)_{i}((X_{k\widetilde{\delta}-\delta}-X^{\delta}_{k\widetilde{\delta}-\delta})^{*}\sigma)_{j}(X^{\delta}_{k\widetilde{\delta}-\delta})(W^{i}_{(k+1)\widetilde{\delta}-\delta}-W^{i}_{k\widetilde{\delta}-\delta})(W^{j}_{(k+1)\widetilde{\delta}-\delta}-W^{j}_{k\widetilde{\delta}-\delta})\\
-&\frac{4}{\gamma}\sum_{i,j=1}^{r}\mu_{\delta}(k\widetilde{\delta}-\delta)(D\varphi\sigma)_{i}((X_{k\widetilde{\delta}-\delta}-X^{\delta}_{k\widetilde{\delta}-\delta})^{*}\sigma)_{j}(X^{\delta}_{k\widetilde{\delta}-\delta})c_{ji}\widetilde{\delta}\\
:=&I_{72121}^{8,k}+I_{72122}^{9,k}.
\end{aligned}
\end{equation}
As in $(\ref{2,2,2,3})$ we can turn $I_{7215}^{k}$ into
$$
I_{7215}^{10,k}:=
\frac{4}{\gamma}\sum_{i,j=1}^{r}\mu_{\delta}(k\widetilde{\delta}-\delta)(D\varphi\sigma)_{i}(X_{k\widetilde{\delta}-\delta})((X_{k\widetilde{\delta}-\delta}-X^{\delta}_{k\widetilde{\delta}-\delta})^{*}\sigma)_{j}(X^{\delta}_{k\widetilde{\delta}-\delta})s(W^{i}_{(k+1)\widetilde{\delta}-\delta}-W^{i}_{k\widetilde{\delta}-\delta})(W^{j}_{(k+1)\widetilde{\delta}-\delta}-W^{j}_{k\widetilde{\delta}-\delta}).
$$
To deal with $I_{723}^{k}$, we write it as
\begin{equation}
\begin{aligned}
&-2\int_{k\widetilde{\delta}}^{(k+1)\widetilde{\delta}}\mu_{\delta}(k\widetilde{\delta}-\delta)(X_{k\widetilde{\delta}-\delta}-X^{\delta}_{k\widetilde{\delta}-\delta})^{*}\left(\int_{k\widetilde{\delta}-\delta}^{u}\bigtriangledown\sigma\sigma(X^{\delta}_{s})dB^{\delta}_{s} \right)dB^{\delta}_{u}\\
&-2\int_{k\widetilde{\delta}}^{(k+1)\widetilde{\delta}}\mu_{\delta}(k\widetilde{\delta}-\delta)(X_{k\widetilde{\delta}-\delta}-X^{\delta}_{k\widetilde{\delta}-\delta})^{*}\left(\int_{k\widetilde{\delta}-\delta}^{u}\bigtriangledown\sigma b(X^{\delta}_{s})\right)dsdB^{\delta}_{u}\\
&-2\int_{k\widetilde{\delta}}^{(k+1)\widetilde{\delta}}\mu_{\delta}(k\widetilde{\delta}-\delta)(X_{k\widetilde{\delta}-\delta}-X^{\delta}_{k\widetilde{\delta}-\delta})^{*}\left(\int_{k\widetilde{\delta}-\delta}^{u}\bigtriangledown\sigma(X^{\delta}_{s})dK^{\delta}_{s}\right)dB^{\delta}_{u}\\
:=&~\sum_{i=1}^{3}I_{723i}^{k},
\end{aligned}
\end{equation}
where
$
(\bigtriangledown\sigma b)_{j}^{i}:=\sum_{\alpha=1}^{d}b^{\alpha}\partial_{x_{\alpha}}\sigma_{j}^{i}.
$
Since $\sigma\in\mathcal{C}_{b}^{2}$, it's easy to see that
$$
E\left[\sup_{s\leq t}|\sum_{k=1}^{S_{\delta}(s)-1}(I_{7232}^{k}+I_{7233}^{k})|^{2}\right]
\lesssim
\left(\sum_{k}E\left[\left(\int_{k\widetilde{\delta}}^{(k+1)\widetilde{\delta}}|\dot{B}^{\delta}_{u}|du\right)^{4}\right]\right)^{\frac{1}{2}}\lesssim n(\delta)^{\frac{3}{2}}\delta^{\frac{1}{2}}.
$$
A similar way as we used in $(\ref{2,2,1,1})$ transforms $I_{7231}^{k}$ to
\begin{equation}
\begin{aligned}
-2&\mu_{\delta}(k\widetilde{\delta}-\delta)\sum_{l=1}^{r}(X_{k\widetilde{\delta}-\delta}-X^{\delta}_{k\widetilde{\delta}-\delta})^{l}\sum_{i,j=1}^{r}\sum_{\alpha=1}^{d}\sigma^{\alpha}_{i}\partial_{x_{\alpha}}\sigma^{l}_{j}(X^{\delta}_{k\widetilde{\delta}-\delta})(W^{i}_{(k+1)\widetilde{\delta}-\delta}-W^{i}_{k\widetilde{\delta}-\delta})(W^{j}_{(k+1)\widetilde{\delta}-\delta}-W^{j}_{k\widetilde{\delta}-\delta})
\\
+2&\mu_{\delta}(k\widetilde{\delta}-\delta)\sum_{l=1}^{r}(X_{k\widetilde{\delta}-\delta}-X^{\delta}_{k\widetilde{\delta}-\delta})^{l}\sum_{i,j=1}^{r}\sum_{\alpha=1}^{d}\sigma^{\alpha}_{i}\partial_{x_{\alpha}}\sigma^{l}_{j}(X^{\delta}_{k\widetilde{\delta}-\delta})c_{ji}\widetilde{\delta}\\
:=~&I_{72311}^{11,k}+I_{72312}^{12,k}.
\end{aligned}
\end{equation}

So far, we have left $I_{6}(t)$, $I_{8}(t)$, $I_{9}(t)$, $I_{10}(t)$ and $I_{2228}^{1,k}$, $I_{2227}^{2,k}$, $I_{2227}^{3,k}$, $I_{7241}^{4,k}$, $I_{7221}^{5,k}$, $I_{72241}^{6,k}$, $I_{72242}^{7,k}$, $I_{72121}^{8,k}$, $I_{72122}^{9,k}$, $I_{7215}^{10,k}$, $I_{72311}^{11,k}$, $I_{72312}^{12,k}$. We are going to bound them.

$Group~1$. $I_{8}$, $I_{7241}^{4,k}$.

Since
\begin{equation}
\begin{aligned}
&~~I_{8}(s)+\sum_{k=1}^{S_{\delta}(s)-1}~I_{7241}^{4,k}\\
=~~&2\sum_{k=1}^{S_{\delta}(s)-1}\int_{k\widetilde{\delta}-\delta}^{(k+1)\widetilde{\delta}-\delta}\mu_{\delta}(u)(X_{u}-X^{\delta}_{u})^{*}\sigma(X_{u})-\mu_{\delta}(k\widetilde{\delta}-\delta)(X_{k\widetilde{\delta}-\delta}-X^{\delta}_{k\widetilde{\delta}-\delta})^{*}\sigma(X_{k\widetilde{\delta}-\delta})dW_{u}\\
+&2\int_{0}^{\widetilde{\delta}-\delta}\mu_{\delta}(u)(X_{u}-X^{\delta}_{u})^{*}\sigma(X_{u})dW_{u}\\
+&2\int_{[s](\widetilde{\delta})-\delta}^{s}\mu_{\delta}(u)(X_{u}-X^{\delta}_{u})^{*}\sigma(X_{u})dW_{u},
\end{aligned}
\end{equation}
BDG's inequality shows
$$
E\left[\sup_{s\leq t}|I_{8}(s)+\sum_{k=1}^{S_{\delta}(s)-1}I_{7241}^{4,k}|^{2}\right]\lesssim n(\delta)^{\frac{1}{2}}\delta^{\frac{1}{2}}.
$$

$Group~2$. $I_{9}$, $I_{72311}^{11,k}$, $I_{72312}^{12,k}$.

Neglecting some trivial terms, we can consider $I_{9}(s)$ as
$$
2\mu_{\delta}(k\widetilde{\delta}-\delta)\sum_{k=1}^{S_{\delta}(s)-1}\sum_{l=1}^{r}(X_{k\widetilde{\delta}-\delta}-X^{\delta}_{k\widetilde{\delta}-\delta})^{l}\sum_{i,j=1}^{r}\sum_{\alpha=1}^{d}\sigma^{\alpha}_{i}\partial_{x_{\alpha}}\sigma^{l}_{j}(X^{\delta}_{k\widetilde{\delta}-\delta})c_{ij}\widetilde{\delta}.
$$
Note
$$
c_{i j}=s_{i j}+\frac{1}{2} \delta_{i j}, \quad i, j=1,2, \ldots, r
$$
and  $\left(s_{i j}\right)$ is a skew-symmetric $r \times r-$matrix, we have
\begin{equation}
\begin{aligned}
&~~I_{9}(s)+\sum_{k=1}^{S_{\delta}(s)-1}I_{72312}^{12,k}\\
=&~~2\mu_{\delta}(k\widetilde{\delta}-\delta)\sum_{k=1}^{S_{\delta}(s)-1}\sum_{l=1}^{r}(X_{k\widetilde{\delta}-\delta}-X^{\delta}_{k\widetilde{\delta}-\delta})^{l}\sum_{i=1}^{r}\sum_{\alpha=1}^{d}\sigma^{\alpha}_{i}\partial_{x_{\alpha}}\sigma^{l}_{i}(X^{\delta}_{k\widetilde{\delta}-\delta})\widetilde{\delta}.
\end{aligned}
\end{equation}
By It\^o's formula, we have
\begin{equation}\label{7,2,3,1,1}
\begin{aligned}
&~~~~\sum_{k=1}^{S_{\delta}(s)-1}I_{72311}^{11,k}\\
=&-2\sum_{k=1}^{S_{\delta}(s)-1}\mu_{\delta}(k\widetilde{\delta}-\delta)\sum_{l=1}^{r}(X_{k\widetilde{\delta}-\delta}-X^{\delta}_{k\widetilde{\delta}-\delta})^{l}\sum_{i=1}^{r}\sum_{\alpha=1}^{d}\sigma^{\alpha}_{i}\partial_{x_{\alpha}}\sigma^{l}_{i}(X^{\delta}_{k\widetilde{\delta}-\delta})\widetilde{\delta}\\
&-2\sum_{k=1}^{S_{\delta}(s)-1}\int_{k\widetilde{\delta}-\delta}^{(k+1)\widetilde{\delta}-\delta}\mu_{\delta}(k\widetilde{\delta}-\delta)\sum_{l=1}^{r}(X_{k\widetilde{\delta}-\delta}-X^{\delta}_{k\widetilde{\delta}-\delta})^{l}
\sum_{i,j=1}^{r}\sum_{\alpha=1}^{d}\sigma^{\alpha}_{i}\partial_{x_{\alpha}}\sigma^{l}_{j}(X^{\delta}_{k\widetilde{\delta}-\delta})(W^{i}_{u}-W^{i}_{k\widetilde{\delta}-\delta})dW^{j}_{u}\\
&-2\sum_{k=1}^{S_{\delta}(s)-1}\int_{k\widetilde{\delta}-\delta}^{(k+1)\widetilde{\delta}-\delta}\mu_{\delta}(k\widetilde{\delta}-\delta)\sum_{l=1}^{r}(X_{k\widetilde{\delta}-\delta}-X^{\delta}_{k\widetilde{\delta}-\delta})^{l}
\sum_{i,j=1}^{r}\sum_{\alpha=1}^{d}\sigma^{\alpha}_{i}\partial_{x_{\alpha}}\sigma^{l}_{j}(X^{\delta}_{k\widetilde{\delta}-\delta})(W^{j}_{u}-W^{j}_{k\widetilde{\delta}-\delta})dW^{i}_{u}.
\end{aligned}
\end{equation}
By BDG's inequality, the last two terms are trivial. However the first one plus
$$
I_{9}(s)+\sum_{k=1}^{S_{\delta}(s)-1}I_{72312}^{12,k}
$$
equals to zero, thus
$$
E\left[\sup_{s\leq t}|I_{9}(s)+\sum_{k=1}^{S_{\delta}(s)-1}(I_{72311}^{11,k}+I_{72312}^{12,k})|^{2}\right]
$$
is trivial.

$Group~3$. $I_{10}$, $I_{7221}^{5,k}$, $I_{72241}^{6,k}$, $I_{72242}^{7,k}$.

Applying the same way in $(\ref{7,2,3,1,1})$ to $I_{7221}^{5,k}$ and $I_{72241}^{6,k}$, we get
$$
-2\mu_{\delta}(k\widetilde{\delta}-\delta)\sum_{i,j=1}^{r}\sigma_{ij}(X_{k\widetilde{\delta}-\delta})\sigma_{ij}(X^{\delta}_{k\widetilde{\delta}-\delta})\widetilde{\delta}
$$
and
$$
2\mu_{\delta}(k\widetilde{\delta}-\delta)\sum_{i,j=1}^{r}\sigma_{ij}^{2}(X^{\delta}_{k\widetilde{\delta}-\delta})\widetilde{\delta}
$$
without some trivial terms.
By virtue of $(\ref{s_{i j}})$,
$$
I_{72242}^{7,k}=-\mu_{\delta}(k\widetilde{\delta}-\delta)\sum_{i,j=1}^{r}\sigma_{ij}^{2}(X^{\delta}_{k\widetilde{\delta}-\delta})\widetilde{\delta}.
$$
Thus we claim
\begin{equation}
\begin{aligned}
&~~~E\left[\sup_{s\leq t}|I_{10}(s)+\sum_{k=1}^{S_{\delta}(s)-1}(I_{7221}^{5,k}+I_{72241}^{6,k}+I_{72242}^{7,k})|^{2}\right]\\
\lesssim&~\sum_{i,j=1}^{r}E\left[\left(\sum_{k=1}^{S_{\delta}(t)-1}\mu_{\delta}(k\widetilde{\delta}-\delta)\left(\sigma_{ij}(X^{\delta}_{k\widetilde{\delta}-\delta})-\sigma_{ij}(X_{k\widetilde{\delta}-\delta})\right)^{2}\widetilde{\delta}\right)^{2}\right]+n(\delta)^{\frac{1}{2}}\delta^{\frac{1}{2}}\\
\lesssim&~\int_{0}^{t}E[m_{\delta}^{2}(u)]du+n(\delta)^{\frac{1}{2}}\delta^{\frac{1}{2}}.
\end{aligned}
\end{equation}

$Group~4$. $I_{2227}^{2,k}$, $I_{2227}^{3,k}$, $I_{72122}^{9,k}$.

It's easy to see that
$$
I_{2227}^{2,k}+I_{72122}^{9,k}
=~-\frac{4}{\gamma}\sum_{i=1}^{r}\mu_{\delta}(k\widetilde{\delta}-\delta)(D\varphi\sigma)_{i}((X_{k\widetilde{\delta}-\delta}-X^{\delta}_{k\widetilde{\delta}-\delta})^{*}\sigma)_{i}(X^{\delta}_{k\widetilde{\delta}-\delta})\widetilde{\delta}.
$$
A similar way as we handle $(\ref{7,2,3,1,1})$, we can turn $I_{2227}^{3,k}$ into
$$
\frac{4}{\gamma}\sum_{i=1}^{r}\mu_{\delta}(k\widetilde{\delta}-\delta)(D\varphi\sigma)_{i}((X_{k\widetilde{\delta}-\delta}-X^{\delta}_{k\widetilde{\delta}-\delta})^{*}\sigma)_{i}(X^{\delta}_{k\widetilde{\delta}-\delta})\widetilde{\delta},
$$
if we ignore some trivial terms. Hence
\begin{equation}
E\left[\sup_{s\leq t}|\sum_{k=1}^{S_{\delta}(s)-1}(I_{2227}^{2,k}+I_{2227}^{3,k}+I_{72122}^{9,k})|^{2}\right]\lesssim o(1).
\end{equation}

$Group~5$. $I_{2228}^{1,k}$, $I_{72121}^{8,k}$.

Note that
\begin{equation}
\begin{aligned}
&~~~\sum_{k=1}^{n}(I_{2228}^{1,k}+I_{72121}^{8,k})\\
=~&\frac{4}{\gamma}\sum_{k=1}^{n}\sum_{i,j=1}^{r}\mu_{\delta}(D\varphi\sigma)_{i}(X^{\delta}_{k\widetilde{\delta}-\delta})\left[((X_{k\widetilde{\delta}-\delta}-X^{\delta}_{k\widetilde{\delta}-\delta})^{*}\sigma)_{j}(X^{\delta}_{k\widetilde{\delta}-\delta})\right.\\
&~~~~\left.-((X_{k\widetilde{\delta}-\delta}-X^{\delta}_{k\widetilde{\delta}-\delta})^{*}\sigma)_{j}(X_{k\widetilde{\delta}-\delta})\right](W^{i}_{(k+1)\widetilde{\delta}-\delta}-W^{i}_{k\widetilde{\delta}-\delta})(W^{j}_{(k+1)\widetilde{\delta}-\delta}-W^{j}_{k\widetilde{\delta}-\delta}).
\end{aligned}
\end{equation}
A similar way as we handle $(\ref{7,2,3,1,1})$, we can turn it into
\begin{equation}
\begin{aligned}
&\frac{4}{\gamma}\sum_{k=1}^{n}\sum_{i=1}^{r}\mu_{\delta}(D\varphi\sigma)_{i}(X^{\delta}_{k\widetilde{\delta}-\delta})\left[((X_{k\widetilde{\delta}-\delta}-X^{\delta}_{k\widetilde{\delta}-\delta})^{*}\sigma)_{i}(X^{\delta}_{k\widetilde{\delta}-\delta})-((X_{k\widetilde{\delta}-\delta}-X^{\delta}_{k\widetilde{\delta}-\delta})^{*}\sigma)_{i}(X_{k\widetilde{\delta}-\delta})\right]\widetilde{\delta}.
\end{aligned}
\end{equation}
Then we get
\begin{equation}
E\left[\sup_{s\leq t}|\sum_{k=1}^{S_{\delta}(s)-1}(I_{2228}^{1,k}+I_{72121}^{8,k})|^{2}\right]\lesssim\int_{0}^{t}E[m_{\delta}^{2}(u)]du+n(\delta)^{\frac{1}{2}}\delta^{\frac{1}{2}}.
\end{equation}

$Group~6$. $I_{6}$, $I_{7215}^{10,k}$.

Assume that
$$
I_{6}^{k}
:=-\frac{4}{\gamma}\sum_{i=1}^{r}\mu_{\delta}(k\widetilde{\delta}-\delta)(D\varphi\sigma)_{i}((X_{k\widetilde{\delta}-\delta}-X^{\delta}_{k\widetilde{\delta}-\delta})^{*}\sigma)_{i}(X_{k\widetilde{\delta}-\delta})\widetilde{\delta}.
$$
Clearly, it holds that
$$
E\left[\sup_{s\leq t}|\sum_{k=1}^{S_{\delta}(s)-1}I_{6}^{k}-I_{6}(s)|^{2}\right]\lesssim n(\delta)^{\frac{1}{2}}\delta^{\frac{1}{2}}.
$$
Moreover, as we do in $(\ref{7,2,3,1,1})$, by Ito's formula,
\begin{equation}
\begin{aligned}
&~~~~\sum_{k=1}^{S_{\delta}(s)-1}~\left(I_{6}^{k}+I_{7215}^{10,k}\right)\\
=&\frac{4}{\gamma}\sum_{k=1}^{S_{\delta}(s)-1}\sum_{i=1}^{r}\mu_{\delta}(k\widetilde{\delta}-\delta)(D\varphi\sigma)_{i}(X_{k\widetilde{\delta}-\delta})\left[((X_{k\widetilde{\delta}-\delta}-X^{\delta}_{k\widetilde{\delta}-\delta})^{*}\sigma)_{i}(X^{\delta}_{k\widetilde{\delta}-\delta})-((X_{k\widetilde{\delta}-\delta}-X^{\delta}_{k\widetilde{\delta}-\delta})^{*}\sigma)_{i}(X_{k\widetilde{\delta}-\delta})\right]\widetilde{\delta}\\
+&\frac{4}{\gamma}\sum_{k=1}^{S_{\delta}(s)-1}\int_{k\widetilde{\delta}-\delta}^{(k+1)\widetilde{\delta}-\delta}\mu_{\delta}(k\widetilde{\delta}-\delta)\sum_{i,j=1}^{r}(D\varphi\sigma)_{i}(X_{k\widetilde{\delta}-\delta})((X_{k\widetilde{\delta}-\delta}-X^{\delta}_{k\widetilde{\delta}-\delta})^{*}\sigma)_{j}(X^{\delta}_{k\widetilde{\delta}-\delta})(W^{i}_{u}-W^{i}_{k\widetilde{\delta}-\delta})dW^{j}_{u}\\
+&\frac{4}{\gamma}\sum_{k=1}^{S_{\delta}(s)-1}\int_{k\widetilde{\delta}-\delta}^{(k+1)\widetilde{\delta}-\delta}\mu_{\delta}(k\widetilde{\delta}-\delta)\sum_{i,j=1}^{r}(D\varphi\sigma)_{i}(X_{k\widetilde{\delta}-\delta})((X_{k\widetilde{\delta}-\delta}-X^{\delta}_{k\widetilde{\delta}-\delta})^{*}\sigma)_{j}(X^{\delta}_{k\widetilde{\delta}-\delta})(W^{j}_{u}-W^{j}_{k\widetilde{\delta}-\delta})dW^{i}_{u}.
\end{aligned}
\end{equation}
BDG's inequality and Lemma $\ref{l8}$ yield
$$
~E\left[\sup_{s\leq t}|\sum_{k=1}^{S_{\delta}(s)-1}\left(I_{6}^{k}+I_{7215}^{10,k}\right)|^{2}\right]\lesssim~\int_{0}^{t}E[m_{\delta}^{2}(u)]du+n(\delta)^{\frac{1}{2}}\delta^{\frac{1}{2}}.
$$
It follows that
$$
~E\left[\sup_{s\leq t}|\sum_{k=1}^{S_{\delta}(s)-1}I_{7215}^{10,k}+I_{6}(s)|^{2}\right]\lesssim~\int_{0}^{t}E[m_{\delta}^{2}(u)]du+n(\delta)^{\frac{1}{2}}\delta^{\frac{1}{2}}.
$$
Since $\varphi\in\mathcal{C}_{b}^{2}$, put estimates of $Group~1-6$ together, we claim
\begin{equation}
\begin{aligned}
	E\left[\sup_{u\leq t}|X_{u}-X^{\delta}_{u}|^4\right]\lesssim~&\int_{0}^{t}E[m_{\delta}^{2}(u)]du+o(1)
	\\
	\lesssim~&\int_{0}^{t}E[\sup_{u\leq s}|X_{u}-X^{\delta}_{u}|^{4}]ds+o(1),
\end{aligned}
\end{equation}
where
$
o(1)=n(\delta)^{\frac{5}{2}}\delta^{\frac{1}{2}}+n(\delta)^{-1}+|(c_{ij})(\widetilde{\delta},\delta)-(c_{ij})|^{2}.
$
We complete the prove.
\end{proof}

$\boldsymbol{Proof~of~Theorem~\ref{t1}}$.
\begin{proof}
Employing Proposition $\ref{p2}$ and  Gronwall's inequality, we get
\begin{equation}\label{theorem}
	E\left[\|X^{\delta,T}-X^{T}\|\right]\leq C\left(n(\delta)^{\frac{5}{2}}\delta^{\frac{1}{2}}+n(\delta)^{-1}+|(c_{ij})(\widetilde{\delta},\delta)-(c_{ij})|^{2}\right).
\end{equation}
That is, $\forall\epsilon>0$,
$$
\lim_{\delta\to 0+}P(\|X^{\delta,T}-X^{T}\|>\epsilon)=\lim_{\delta\to 0+}\frac{E\left[\|X^{\delta,T}-X^{T}\|\right]}{\epsilon}=0.
$$
Furthermore, by the boundness of solution on $D$ and Dominated Convergence Theorem, we have
$$
\lim_{\delta\to 0+}E\left[\|X^{\delta,T}-X^{T}\|^p\right]=0
$$
for any $p>0$.
\end{proof}

\section{Example}

~~~~~According to Theorem $\ref{t1}$, we can get strong convergence of some usual approximation.

\begin{example}$($mollifiers$)$.
Let $\rho$ be a non-negative $C^{\infty}-$ function whose support is contained in $[0,1]$, satisfying $\int_{0}^{1} \rho(s) d s=1$. we define
$$
\rho_{\delta}(s)=\frac{1}{\delta} \rho\left(\frac{s}{\delta}\right) \quad \text { for } \quad \delta>0
$$
and
$$
G^{\delta,i}(t, w)=\int_{0}^{\infty} w^{i}(s) \rho_{\delta}(s-t) d s=\int_{0}^{\delta} w^{i}(s+t) \rho_{\delta}(s) d s
$$
for $i=1,2, \ldots, r$.
\end{example}
\begin{proof}
It's easy to see that  $\{G^{\delta}(t, f)\}_{\delta>0}$ is a class of $\mathcal{B}(R)\times\mathbb{C}[0,\infty)-$measurable map.

$(i)$ for any $t\in[0,\infty)$, we have
$$
\dot{G}^{\delta,i}(t, w)=-\frac{1}{\delta} \int_{0}^{1} w^{i}(t+\delta \xi) \rho^{\prime}(\xi) d \xi\in C[0,\infty).
$$

$(ii)$ is obvious.

$(iii)$
\begin{equation}
\begin{aligned}
	G^{\delta,i}(t+k\delta, w)=&\int_{0}^{\delta} w^{i}(s+t+k\delta) \rho_{\delta}(s) d s\\
	=&\int_{0}^{\delta} \left(w^{i}(s+t+k\delta)-w^{i}(k\delta)\right) \rho_{\delta}(s) d s+w^{i}(k\delta)\\
	=&~G^{\delta,i}(t,\theta_{k\delta}w)+w^{i}(k\delta).
\end{aligned}
\end{equation}

$(iv)$
$$
E\left[G^{\delta,i}(0, w)\right]=\int_{0}^{\delta} E\left[W^{i}_{s}\right] \rho_{\delta}(s) d s=0.
$$

$(v)$
\begin{equation}
\begin{aligned}
	&\delta^{m-1}E\left[\int_{0}^{\delta}|\dot{G}^{\delta,i}(s, w)|^{2m}ds \right]
	=\frac{1}{\delta^{m+1}}E\left[\int_{0}^{\delta}|\int_{0}^{1} w^{i}(s+\delta \xi) \rho^{\prime}(\xi) d \xi|^{2m}ds\right]\\
	=&E\left[\int_{0}^{1}|\int_{0}^{1}\frac{w^{i}(\delta s+\delta\xi)}{\sqrt{\delta}}\rho^{\prime}(\xi) d \xi|^{2m}ds\right]
	=E\left[\int_{0}^{1}|\int_{0}^{1}w^{i}(s+\xi)\rho^{\prime}(\xi) d \xi|^{2m}ds\right]<\infty.
\end{aligned}
\end{equation}

$(vi)$ By H\"older's inequality, we have
\begin{equation}
\begin{aligned}
&E\left[\left|G^{\delta,i}(0,f(\ast))\right|^{2p}\right]=E\left[\left|\int_{0}^{\delta} f(s) \rho_{\delta}(s) d s\right|^{2p}\right]
\leq\int_{0}^{\delta} E\left[|f(s)|^{2p}\right]ds \left(\int_{0}^{\delta}\rho_{\delta}^{\frac{2p}{2p-1}}(s)ds\right)^{2p-1}\\
=~&\delta
\left(\int_{0}^{\delta}\rho_{\delta}^{\frac{2p}{2p-1}}(s)ds\right)^{2p-1}E\left[|f(\xi)|^{2p}\right]
=
\left(\int_{0}^{1}\rho^{\frac{2p}{2p-1}}(s)ds\right)^{2p-1}E\left[|f(\xi)|^{2p}\right],
\end{aligned}
\end{equation}
for some $\xi\in[0,\delta]$. Finally, we have $s_{ij}(\delta,\delta)=0$. Thus we say $\{G^{\delta}(t, f)\}_{\delta>0}$ is an approximation of Brownian Motion.
\end{proof}

Now, let $C^{1}_{0,1}$ be the space of continuously differentiable functions $f$ on $[0,1]$ such
that $f(0)=0, f(1)=1 .$ Let $f'=\frac{d}{d t} f$ and $\Delta_{k+1} w=w(k \delta+\delta)-w(k\delta)$.

\begin{example} Select $f^{i} \in C^{1}_{0,1}, i=1,2, \ldots, r$, and set
$$
G^{\delta,i}(t, w)=w^{i}_{t_{\delta}}+f^{i}\left(\frac{t-t_{\delta}}{\delta}\right)(w^{i}_{\overline{t}_{\delta}}-w^{i}_{t_{\delta}}) .
$$
\end{example}
\begin{proof}
$(i)-(iii)$ is obvious. Since $G^{\delta,i}(0, w)\equiv0$, $(iv)$ and $(vi)$ are also easy to get. As for $(iv)$, we have
\begin{equation}
\begin{aligned}
	\delta^{m-1}E\left[\int_{0}^{\delta}|\dot{G}^{\delta,i}(s, w)|^{2m}ds \right]
	=&~\frac{1}{\delta^{m+1}}\int_{0}^{\delta}\left|\dot{f}^{i}\left(\frac{s}{\delta}\right)\right|^{2m}ds~E\left[\left|W^{i}_{\delta}-W^{i}_{0}\right|^{2m}\right]\\
	\lesssim&~\int_{0}^{1}|\dot{f}^{i}(s)|^{2m}ds<\infty.
\end{aligned}
\end{equation}
Since $s_{ij}(\delta,\delta)=0$, Assumption $\ref{a1}$ is satisfied.
\end{proof}

\begin{example} $($McShane$)$. If $r=2$ and choose $f^{i} \in C^{1}_{0,1}$, $i=1,2$. Let
\begin{equation}
G^{\delta,i}(t, w)=
\left\{
\begin{array}{cc}
w^{i}(t_{\delta})+f^{i}\left(\frac{t-t_{\delta}}{\delta}\right)(w^{i}_{\overline{t}_{\delta}}-w^{i}_{t_{\delta}}), & (w^{1}_{\overline{t}_{\delta}}-w^{1}_{t_{\delta}})(w^{2}_{\overline{t}_{\delta}}-w^{2}_{t_{\delta}}) \geq 0,\\
w^{i}(t_{\delta})+f^{3-i}\left(\frac{t-t_{\delta}}{\delta}\right)(w^{i}_{\overline{t}_{\delta}}-w^{i}_{t_{\delta}}), & (w^{1}_{\overline{t}_{\delta}}-w^{1}_{t_{\delta}})(w^{2}_{\overline{t}_{\delta}}-w^{2}_{t_{\delta}})<0.
\end{array}
\right.
\end{equation}
\end{example}
\begin{proof}
A similar argument shows $(i)-(vi)$ are satisfied, now we say that Assumption $\ref{a1}$ is also true. In fact, using integration by parts, we have
$$
\frac{1}{2}\int_{0}^{\delta}[G^{\delta,1}(s, w)\dot{G}^{\delta,2}(s, w)-G^{\delta,2}(s, w)\dot{G}^{\delta,1}(s, w)]ds
=\frac{1-2\int_{0}^{1}\dot{f}^{1}(s)f^{2}(s)ds}{2}|W^{1}_{\delta}||W^{2}_{\delta}|.
$$
\end{proof}

\end{document}